\documentclass[11pt]{article}

\usepackage{lscape,tikz}

\definecolor{darkblue}{rgb}{0.2,0.2,0.71}
\usepackage{todonotes}
\usepackage{epsfig,cancel,ulem, amsthm,amssymb}
\usepackage{color,soul}
\usetikzlibrary{decorations.markings,arrows}

\usepackage{wrapfig}
\usepackage{color}
\usepackage{graphicx,framed}
\definecolor{shadecolor}{rgb}{0.95, 0.95, 0.86}
\definecolor{darkgreen}{rgb}{0.2, 0.5,  0}

\newcommand{\rk}{\mathrm{rank}}

\def\&{\vspace{-5pt}&}

\def\Tr{ {\rm Tr}}

\def \eqref#1{(\ref{#1})}
\def \& {&\hspace{-10pt}}

\newcommand{\al}{\alpha}

\newcommand{\bt}{\beta}

\newcommand{\ga}{\gamma}

\renewcommand{\O}{\Omega}

\newcommand{\br}{{\mathbb R}}

\newcommand{\nn}{\nonumber}

\newcommand{\g}{\mathfrak{g}}

\newcommand{\fn}{\mathfrak{n}}

\newcommand{\Y}{\Lambda}
 \newcommand{\ew}{\mathrm{W}}
 \newcommand{\exx}{\mathrm{E}}
  \newcommand{\A}{\mathcal{A}}
    \newcommand{\B}{\mathcal{B}}
    \newcommand{\V}{\mathcal{V}}
     \newcommand{\lop}[1]{\mathfrak{L}(#1)}
      \newcommand{\Lie}{\mathfrak{L}}
\newcommand{\bil}[2]{{\langle #1 | #2\rangle}}
\newcommand{\bin}[2]{{{#1}\choose{#2} }}

\newcommand{\nneg}{\mathfrak{n}}

 \newcommand{\lcan}{\mathcal{L}^{\!c}}
  \newcommand{\bneg}{\mathfrak{b}}
  
 \newcommand{\px}{\partial_x}
     \newcommand{\W}{\mathcal{W}}

\newcommand{\f}{\mathcal{F}}
\newcommand{\I}{\mathcal{I}}
\newcommand{\Q}{\mathcal{Q}}
\newcommand{\N}{\mathcal{N}}

\newtheorem{theorem}{Theorem}[section]
\newtheorem{example}[theorem]{Example}
\newtheorem{exercise}[theorem]{Exercise}

\newtheorem{lemma}[theorem]{Lemma}
\newtheorem{remark}[theorem]{Remark}

\newtheorem{proposition}[theorem]{Proposition}
\newtheorem{corollary}[theorem]{Corollary}
\newtheorem{definition}[theorem]{Definition}

\def\V {\mathcal V}

\def\bt{\begin{theorem}}
\def\et{\end{theorem}}
\def\bc{\begin{corollary}}
\def\ec{\end{corollary}}
\def\bx{\begin{example}\small}
\def\ex{\end{example}}
\def\bxr{\begin{exercise}\small}
\def\exr{\end{exercise}}
\def\bl{\begin{lemma}}
\def\el{\end{lemma}}
\def\bd{\begin{definition}}
\def\ed{\end{definition}}
\def\bp{\begin{proposition}}
\def\ep{\end{proposition}}

\def\br{\begin{remark}}
\def\er{\end{remark}}

\def\be{\begin{equation}}
\def\ee{\end{equation}}

\def\beq{\begin{equation}}
\def\eeq{\end{equation}}
\def\&{\hspace{-15pt}&}
\def\bea{\begin{eqnarray}}
\def\eea{\end{eqnarray}}

\def\L{\mathcal L}

\def\H{{\cal H}}

\def\1{{\bf 1}}

\newcommand{\h}{\mathfrak{h}}

\newcommand{\ad}{\mathrm{ad}}

\makeatletter
\@addtoreset{equation}{section}
\makeatother

\begin{document}
\title{Algebraic classical $W$-algebras and Frobenius manifolds}
\author{Yassir Ibrahim Dinar}
\date{}
\maketitle
{\begin{center}{Dedicated to the memory of Boris} Dubrovin\end{center}}
\begin{abstract}
We consider Drinfeld-Sokolov  bihamiltonian structure associated to a distinguished nilpotent elements of semisimple type and the space of common equilibrium points defined by its leading term. On this space, we construct a local bihamiltonian structure which forms an exact Poisson pencil, defines an algebraic classical $W$-algebra,  admits a dispersionless limit, and its leading term defines  an algebraic Frobenius manifold.  This leads to a uniform construction of algebraic Frobenius manifolds corresponding  to regular cuspidal conjugacy classes in irreducible Weyl groups.
\end{abstract}

\tableofcontents

\setcounter{equation}{0}
\setcounter{theorem}{0}
\section{Introduction}\label{intro-s}
\setcounter{equation}{0}
\setcounter{theorem}{0}

Frobenius manifold is a marvelous geometric realization introduced by Boris Dubrovin for undetermined  partial differential equations known as Witten-Dijkgraaf-Verlinde-Verlinde (WDVV) equations which describe the  module space of two dimensional topological field theory. Remarkably, Frobenius manifolds are also  recognized  in many other fields in mathematics like  invariant theory, quantum cohomology, integrable systems  and singularity theory \cite{DuRev}. Briefly,
a Frobenius manifold is a manifold  with a smooth structure of Frobenius algebra on the tangent space   with certain compatibility conditions. By Frobenius algebra, we mean a commutative associative algebra with unity  and an invariant nondegenerate symmetric bilinear form.

Let $M$ be a Frobenius manifold. Then,  we  require the  bilinear form $(.,.)$ to be flat and the unity vector field $e$ is constant with respect to it. Let $(t^1,\ldots,t^r)$  be flat coordinates for $(.,.)$ where  $e={\partial_{t^{r}}}$. Then  the compatibility conditions imply  that there exists a function $\mathbb{F}(t^1,\ldots,t^r)$ such that
\begin{equation}\label{flat metric} \eta_{ij}:=(\partial_{t^i},\partial_{t^j})=  \partial_{t^r}
\partial_{t^i}
\partial_{t^j} \mathbb{F}(t)\end{equation}
and the structure constants of the Frobenius algebra are given by
\begin{equation}  C_{ij}^k(t):=\sum_p \eta^{kp}  \partial_{t^p}\partial_{t^i}\partial_{t^j} \mathbb{F}(t)\end{equation}
where the matrix $\eta^{ij}$ is  the inverse of the matrix $\eta_{ij}$. Associativity  in $T_tM$
implies that  $\mathbb{F}(t)$ satisfies  WDVV equations \cite{wdvv}:
\begin{equation} \label{frob}
 \sum_{k,p}\partial_{t^i}
\partial_{t^j}
\partial_{t^k} \mathbb{F}(t)~ \eta^{kp} ~\partial_{t^p}
\partial_{t^q}
\partial_{t^s} \mathbb{F}(t) = \sum_{k,p}\partial_{t^s}
\partial_{t^j}
\partial_{t^k} \mathbb{F}(t) ~\eta^{kp}~\partial_{t^p}
\partial_{t^q}
\partial_{t^i} \mathbb{F}(t),
  \end{equation}
 for all $i,j,q$ and $s$. In this article,  we consider Frobenius manifolds where the quasihomogeneity  condition for $\mathbb{F}(t)$ can take the form
\begin{equation}\label{quasihomog}
\sum_{i=1}^r d_i t^i \partial_{t^i} \mathbb{F}(t) = \left(3-d \right) \mathbb{F}(t);~~~~d_{r}=1.
\end{equation}
The vector field $E= \sum_{i=1}^n d_i t^i \partial_{t_i}$ is known as  Euler vector field and it  defines  the degrees $d_i$ and  the charge $d$ of  $M$. A Frobenius manifold is called algebraic if $\mathbb{F}(t)$ is an algebraic function and it is called semisimple if $T_t M$ is a semisimple algebra for generic point $t$.

This work is related to a conjecture due to Dubrovin which states that semisimple  irreducible algebraic Frobenius manifolds
with positive degrees correspond to  primitive (quasi-Coxeter) conjugacy classes of irreducible finite Coxeter groups \cite{DMD}.   A primitive conjugacy class in  a Coxeter group is  a conjugacy class  which has no representative in a proper Coxeter subgroup (see \cite{CarClassif} for the classification).  Coxeter conjugacy class is an example of  a primitive conjugacy class which exists in any Coxeter group (it is formed by the product of simple reflections in the case of standard reflection representation).  The  conjecture arises from studying the algebraic solutions to associated equations of isomonodromic deformation of an algebraic Frobenius manifold \cite{DMD}. It  leads to a primitive conjugacy class in a Coxeter group by considering the  classification of  finite orbits of the braid group action on tuple of reflections \cite{STF}. A stage to verify the conjecture is to show the existence of these algebraic Frobenius manifolds.

Under the conjecture, it is known that polynomial Frobenius manifolds correspond to  Coxeter conjugacy classes.  Dubrovin   constructed these polynomial Frobenius structures on orbit spaces of the standard reflection representations of Coxeter groups \cite{DCG}. Their isomonodromic deformations lead to Coxeter conjugacy classes \cite{DMD} and C. Hertling \cite{HER} proved (as also conjectured by Dubrovin) that they exhaust the set of all possible polynomial structures up to an equivalence. This classification and other examples  reveal a relation between  orders and eigenvalues of the  conjugacy classes, and charges and degrees of algebraic Frobenius manifolds. More precisely, if the order of a primitive conjugacy class is $\eta_r+1$ and the eigenvalues are $\exp {2\eta_i \pi \mathbf i\over \eta_r+1}$ then the charge of the corresponding Frobenius structure is $\eta_r-1\over \eta_r+1$ and the degrees are $\eta_i+1\over \eta_r+1$.  We depend on this  relation in constructing   algebraic Frobenius structures.

One of the main methods  to obtain examples of Frobenius manifolds exists within the theory of flat pencils of metrics (equivalently, nondegenerate compatible Poisson brackets of hydordynamics type). Besides, the leading terms of  certain type of local compatible Poisson brackets (a local bihamiltonian structure) which admit(s) a dispersionless limit form a flat pencil of metric \cite{DFP}.

One of the  main ideas to find algebraic Frobenius structures is to restrict ourselves to irreducible Weyl groups, i.e.,  crystallographic Coxeter groups, and to consider the associated simple Lie algebras. Then, under the notion of opposite Cartan subalgebra, regular primitive conjugacy classes correspond to certain   nilpotent  orbits of semisimple type. On the other hand, we can obtain compatible local Poisson brackets for any nilpotent orbit using Drinfeld-Sokolov reduction. These Poisson brackets form an exact Poisson pencil and one of them  is  (or satisfies identities leading to) a classical $W$-algebra. However, they  admit a dispersionless limit only when the nilpotent orbit is regular (which corresponds  to Coxeter conjugacy class). In this article, we will work with a larger type of  conjugacy classes called cuspidal. A cuspidal conjugacy class  has no representative in a Coxeter subgroup of smaller rank.  Regular cuspidal conjugacy classes  correspond to what is called distinguished nilpotent orbits of semisimple type \cite{DelFeher}, \cite{Elash}. In other words, we get certain Drinfeld-Sokolov bihamiltonian structures associated  to regular cuspidal conjugacy classes in irreducible Weyl groups.

Examples of  Frobenius manifolds constructed using Drinfeld-Sokolov bihamiltonian structure  can be traced back to the work of I. Krichever \cite{Krich}. In our terminologies, he treated the case of  Coxeter conjugacy classes in Weyl groups of type $A_r$ (here, classical $W$-algebras  are known as  second Gelfand-Dickey brackets).  In \cite{mypaper1}, we gave a generalization to all Coxeter conjugacy classes in  Weyl groups which, as expected,  lead to the  polynomial  Frobenius manifolds.

For  regular  primitive non-Coxeter conjugacy classes, we always get  algebraic non-polynomial Frobenius structures.  Pavlyk obtained  the first example which is related to the Weyl group of type $D_4$ \cite{PAV}. In \cite{mypaper}, we got another example  working with Weyl group of type $F_4$.  We added another 3 by giving a uniform construction related to  certain conjugacy classes in Weyl groups of type $E_r$, $r=6,7,8$ \cite{mypaper3}. In all these cases, we have to perform  Dirac reduction for  the Drinfeld-Sokolov bihamiltonian structure  to   a subspace to get a bihamiltonian structure admitting a dispersionless limit. In this article, we give a  slightly better interpertation for this subspace which leads to a uniform construction of algebraic Frobenius structures for all regular cuspidal conjugacy classes. Precisely, we will prove the following theorem.

\bt \label{main thm}
Let $\g$ be a complex simple Lie algebra of rank $r$. Fix a regular cuspidal conjugacy class  $[w]$  in the Weyl group $\W (\g)$ of $\g$. Assume the order of representatives in $[w]$ is  $\eta_r+1$ and eigenvalues are $\epsilon^{\eta_i}$, $i=1,\ldots,r$, where $\epsilon$ is a primitive $(\eta_r+1)$th root of unity. Let $\mathcal{O}_{L_1}$ be the distinguished  nilpotent orbit of semisimple type  associated to $[w]$  under the notion of opposite Cartan subalgebra. Consider  the finite bihamiltonian structure formed by the leading term of Drinfled-Sokolov bihamiltonain structure associated to a representative $L_1$ of $\mathcal{O}_{L_1}$. Then its space of common equilibrium points   acquires an algebraic Frobenius manifold structure with charge $\frac{\eta_r-1}{\eta_r+1}$ and degrees $\frac{\eta_i+1}{\eta_r+1}$. This structure depends only on the conjugacy class.
\et

We explain in some details the major steps to prove  theorem \ref{main thm} which lead us to a construction of algebraic classical $W$-algebras admitting a dispersionless limit. Let $\g$ be a complex simple Lie algebra of rank $r$ with the Lie bracket $[\cdot,\cdot]$. Define the adjoint representation $\ad: \g \to \textrm{End}(\g)$ by $\ad_{g_1}(g_2):=[g_1,g_2]$. For $g\in \g$, let $\g^g$ denotes the centralizer of $g$ in $\g$, i.e., $\g^g:=\ker \ad_g$. Fix a distinguished nilpotent element $L_1$ of semisimple type (more details are given is section 3 below). Then, using Jacobson-Morozov theorem, we fix a nilpotent element $f$ and a semisimple element $h$ such that $A:=\{L_1,h,f\}\subseteq \g$ is  a $sl_2$-triple with relations \begin{equation}\label{sl2:relation1}
[h,L_1]= L_1,\quad [h,f]=-
f,\quad [L_1,f]= 2 h.
\end{equation}
We normalize the Killing form on $\g$ to get an invariant bilinear form $\bil . . $ such that $\bil {L_1} f=1$.

Let  $\eta_r$ denotes the maximal eigenvalue of $\ad_h$ acting on $\g$. By definition, we  can (and we will) fix  an element $K_1$ for $L_1$ such that  $\ad_h K_1=-\eta_r K_1$ and $h':=L_1 +K_1$ is a regular semisimple element. Thus, $\h':=\ker \ad_{h'}$ is a Cartan subalgebra known as opposite Cartan subalgebra.   The  adjoint group element $w:=\exp {2\pi \mathbf i\over \eta_r+1} \ad_h$  acts on $\h'$ as a representative of regular cuspidal conjugacy class of order $\eta_r+1$ in the underline  Weyl group $\W(\g)$ (see  \cite{Elash},  and the appendix of \cite{DelFeher}).

Let $\eta_1\leq \ldots\leq \eta_r$  be  natural numbers such that $\epsilon^{\eta_i}$ are eigenvalues of $w$ where $\epsilon$ is $(\eta_r+1)$th root of unity.  Let $n=\dim \g^f$, then using representation theory of $sl_2$-subalgebras, there exist natural numbers $\eta_{r+1},\ldots,\eta_n$ such that the eigenvalues of $\ad_h$ on $\g^f$ are $-\eta_i$, $i=1,\ldots,n$. We list all distinguished nilpotent elements of semisimple type in simple Lie algebras and the numbers $\eta_i$ in table \ref{we} below.

 We fix Slodowy slice $Q:=L_1+ \g^f$ as a transverse subspace to the orbit space of $L_1$ at $L_1$. Let  $\lop {\g^f}$ denotes the space of smooth functions from the circle to $\g^f$. The affine  loop space $\Q:=L_1+\mathfrak L(\g^f)$ carries   compatible local Poisson structures (Drinfeld-Sokolov bihamiltonian structure formed by) $\mathbb B_2^\Q$ and $\mathbb B_1^\Q$, where $\mathbb B_2^\Q$ is a classical $W$-algebra \cite{BalFeh1},\cite{fehercomp}
 and $\mathbb B_1^\Q$ is related to a 2-cocycle on $\g$ provided by $K_1$. They  depend only on the adjoint orbit  of $L_1$ and they can be obtained equivalently by using Drinfeld-Sokolov reduction, bihamiltonian reduction and Dirac reduction \cite{mypaper4}. Note that performing any of these reductions, we need to fix a transverse subspace. However, taking a different subspace than $Q$ will lead to isomorphic bihamiltonian structures. As it is already  known by experts,  we will prove in proposition \ref{exactness} that  $\mathbb B_2^\Q$ and $\mathbb B_1^\Q$ form an exact  Poisson pencil.

We identify Slodowy slice $Q$  with the subspace of constant loops of $\Q$.  We can (and will) fix   coordinates $(z^1,\ldots, z^{n})$ for $Q$ such that
\be
Q=L_1+ \sum z^i \gamma_i,~ \gamma_i\in \g^f, \ad_h \gamma_i=- \eta_i \gamma_i,~ i=1,\ldots, n
\ee
where $\gamma_1=f$ and for $q\in Q$,  $z^1=\bil {L_1} q$. Then the leading terms of $\mathbb B_m^\Q,~ m=1,2$,  can be written as follows
\begin{eqnarray}
  \{z^i(x),z^j(y)\}^{[-1]}_m &=& F^{ij}_m(z(x))\delta(x-y),\\\nonumber
  \{z^i(x),z^j(y)\}^{[0]}_m &=& \Omega^{ij}_m(z(x)) \delta' (x-y)+ \sum_k\Gamma_{k,m}^{ij}(z(x)) z_x^k \delta (x-y).
\end{eqnarray}
Such a local Poisson bracket admits a dispersionless limit iff $F^{ij}_m=0$. In general, $F^{ij}_2(z)$ and $F^{ij}_1(z)$ define compatible Poisson structures $B_2^Q$ and $B_1^Q$, respectively, on $Q$. Moreover, $B_2^Q$ can be identified with the transverse Poisson structure of Lie-Poisson structure on $\g$ \cite{mypaper4}.  We assign $\deg z^i=\eta_i+1$. Then after certain normalization, we will prove the following theorem

\bt \label{nice coordinates} There exists a quasihomogenous  change of coordinates on $Q$ in the form
 \be
 t^i=
 \left\{
  \begin{array}{ll}
z^1,&\hbox{i=1,} \\
  z^i+\mathrm{non~ linear~ terms},&\hbox{i=2,\ldots,r,} \\
z^i,&\hbox{i=r+1,\ldots,n.}
  \end{array}
\right.
\ee
 such that
 \begin{enumerate}
     \item $\deg t^i=\deg z^i=\eta_i+1$
     \item $t^1,\ldots,t^r$ form a complete set of  Casimirs of $B_1^Q$ and they are in involution with respect to $B_2^Q$.
 \end{enumerate}
 \et

  We will keep the notations $(t^1,\ldots,t^n)$ for the coordinates obtained in the last theorem (except in section \ref{Alg Frob mani}, we can and will assume they are flat coordinates of the resulted Frobenius structure).

  We are interested in the space of common equilibrium points $N$ of the bihamiltonian structure formed by $B_2^Q$ and $B_1^Q$.
Combining results from  \cite{bolv1} and \cite{mypaper5}, we explain in theorem  \ref{rank of B1} that the argument shift method leads to a completely integrable system for $B_2^Q$ and
\be
 N =\{q\in Q: \ker B_1^Q(q)=\ker B_2^Q(q)\}.
\ee
 Using Chevalley's theorem, we fix homogeneous set of generators   $P_1,\ldots, P_r$  of the ring of invariant polynomials of $\g$ under the adjoint group action. Let $\overline P_i^0$ denotes the restriction of $P_i$ to $Q$. We can choose  $P_1,\ldots, P_r$ such that the following theorem is valid.  Here, we assume $L_1$ is of type $Z_r(a_s)$ where $Z_r$ is the type of $\g$.

 \bt \label{des of N}
 The space of common equilibrium points $N$  is given by
 \begin{eqnarray}
  N&=&\{t: F^{i\beta}_2(t)=0;~ i=1,\ldots,r,~\beta=r+1,\ldots,n\},\label{defEqN}\\
  &=& \{t: \partial_{t^\beta} \overline P_j^0(t)=0; j=r-s+1,\ldots,r,~~\beta=r+1,\ldots,n\}.\label{defEqN1}
 \end{eqnarray}
 Moreover, $(t^1,\ldots,t^r)$ provide local coordinates around generic points of  $N$. In addition, Dirac reduction of the Poisson pencil  $B_\lambda^Q:=B_2^Q+\lambda B_1^Q$  to $N$ is well defined  and leads to the  trivial  Poisson bracket.
 \et

 Then we construct compatible local Poisson brackets on the loop space $\N=\lop N$.

\bt \label{classical in N}
The Dirac reduction of the Poisson pencil  $\mathbb B_\lambda^\Q:=\mathbb B_2^\Q+\lambda \mathbb B_1^\Q$ to $\N$ is well defined and leads to compatible local Poisson brackets $\{.,.\}_\alpha^\N$, $\alpha=1,2$ which admit  a dispersionless limit and form an exact Poisson pencil. Moreover, $\{.,.\}_2^\N$ is an algebraic classical $W$-algebra.
\et

Let us emphasis that theorem \ref{classical in N} implies that the leading terms of the local Poisson brackets on $\N$  are Poisson brackets of hydrodynamic types, i.e.,
\be\label{reducedPB notations1}
  \{t^u(x),t^v(y)\}^{[0]}_\alpha = \O^{uv}_\alpha(t(x)) \delta' (x-y)+ \Gamma_{\alpha k}^{uv}(t(x)) t_x^k \delta (x-y),~ u,v=1,\ldots r, ~\alpha=1,2.\ee
where $t^k, k>r$ are  solutions  of the polynomial equations \eqref{defEqN} defining $N$.

One of the important steps on the construction is to prove that the matrices $\O^{uv}_\alpha(t)$ are nondegenerate. We will show  that this condition follows from the fact that the restriction of the Killing form   on $\g$ to the  Cartan subalgebra $\h'$ is nondegenerate (see proposition \ref{nondeg} below).

Then we will  prove the following.

\bt \label{my thm}
  The two metrics $\O^{uv}_{1}$ and $   \O^{uv}_{2}$  form a flat pencil of metrics on  $N$ which is regular quasihomogeneous of degree $d={\eta_r-1\over \eta_r+1}$.
\et

 In the end, using theorem \ref{dub flat pencil}  due to Dubrovin, we get the proof of theorem \ref{main thm}.


We organize the article as follows. In section 2, we fix notations and terminologies within the theory of local Poisson brackets, flat pencils of metrics and Frobenius manifolds. We review the classification of distinguished nilpotent orbits of semisimple type in simple Lie algebras in section 3 and we will drive some algebraic properties associated to them.   In section 4, we will study  the space $N$ of common equilibrium points  and prove theorems \ref{nice coordinates} and \ref{des of N}.  We review the Drinfeld-Sokolov reduction in section 5 and prove theorem \ref{classical in N}.  In section 6, we will prove theorem \ref{main thm} and we give examples. The notations given in the introduction are in agreement  with  the flow of the article.

\section{Preliminaries}

In this section,  we  recall  relations between local bihamiltonian structures, flat pencils of metrics and Frobenius manifolds. We also review the notion of Dirac reduction for local Poisson brackets.

\subsection{Contravariant metrics and local Poisson brackets} \label{first-sec}

Let $M$ be a smooth manifold of dimension $n$ and  fix  local coordinates $(u^1,\ldots, u^n)$ on $M$. Here, and in what follows, summation with respect to repeated
upper and lower indices is assumed, i.e., We will adopt Einstein summation convention.

\bd A symmetric bilinear form $(. ,. )$ on $T^*M$ is called a contravariant
metric if it is invertible on an open dense subset $M_0 \subseteq M$. We define the contravariant Levi-Civita connection or Christoffel symbols $\Gamma^{ij}_k$  for a contravariant
metric $(. ,. )$ by
\beq
\Gamma^{ij}_k:=-g^{is} \Gamma_{sk}^j
\eeq
where $\Gamma_{sk}^j$ are the  Christoffel symbols of the metric $<. ,. >$ defined on $TM_0$ by the inverse of the matrix $\O^{ij}(u)=(du^i, du^j)$.
We say the metric $(.,.)$ is flat if  $<. ,. >$ is flat.
\ed

Let $(. ,. )$ be a contraviariant metric on $M$ and set $\O^{ij}(u)=(du^i, du^j)$.  Then  we will use $\Omega^{ij}$ to refer to both the metric and the entries defined by the metric. In particular,  Lie derivative of $(. ,. )$ along a vector field $X$  will be written $\Lie_X \Omega^{ij}$ while  $X\Omega^{ij}$   means the vector field $X$ acting on the entry $\Omega^{ij}$.


The loop space $\lop M$ of $M$ is the space of smooth functions from the circle to $M$. A local Poisson bracket $\{.,.\}$ is  a certain bracket on the space of local  functional on $\lop M$ \cite{DZ}. We can write $\{.,.\}$  as a finite summation   of the form
\begin{eqnarray} \label{genLocPoissBra}\label{genLocBraGen}\{u^i(x),u^j(y)\}&=&
\sum_{k=-1}^\infty  \{u^i(x),u^j(y)\}^{[k]}\\\nonumber
\{u^i(x),u^j(y)\}^{[k]}&=&\sum_{l=0}^{k+1} A_{k,l}^{i,j}(u(x))
\delta^{(k-l+1)}(x-y),
 \end{eqnarray}
where $A_{k,l}^{i,j}(u(x))$ are quasihomogeneous polynomials in $\partial_x^m
u^i(x)$ of degree $l$ when we  assign degree
$\partial_x^m u^i(x)$  equals $m$, and $\delta(x-y)$ is the Dirac  delta function defined by
\be \int_{S^1} f(y) \delta(x-y) dy=f(x).\ee

\bd \cite{fehercomp}
A  local Poisson bracket $\{.,.\}$ in the form  \eqref{genLocBraGen} is called a classical $W$-algebra if there exist  local coordinates $(z^1,\ldots,z^n)$ such that
\begin{eqnarray}\label{walgebra2}
\{z^1(x),z^1(y)\}&=&  c \delta^{'''}(x-y) +2 z^1(x) \delta'(x-y)+ z^1_x\delta(x-y), \\\nonumber
\{z^1(x),z^i(y)\} &=& (\eta_i+1) z^i(x) \delta'(x-y)+ \eta_i z^i_x \delta(x-y),
\end{eqnarray}
for nonzero constant $c$.

\ed

Let us fix a local Poisson bracket  $\{.,.\}$ on $\lop M$. The first terms can be written
as follows
\begin{eqnarray}\label{loc poiss}
  \{u^i(x),u^j(y)\}^{[-1]} &=& F^{ij}(u(x))\delta(x-y), \\\nonumber
  \{u^i(x),u^j(y)\}^{[0]} &=& \O^{ij}(u(x)) \delta' (x-y)+ \Gamma_k^{ij}(u(x)) u_x^k \delta (x-y),\\\nonumber
  \{u^i(x),u^j(x)\}^{[k]}& = & S^{ij}_k(u(x)) \delta^{k+1}(x-y)+\ldots,~  ~k>0.
\end{eqnarray}
   Note that $M$ can be defined as the subspace  of constant loops of $\lop M$. Then  $\O^{ij}(u)$, $F^{ij}(u)$, $S^{ij}_k(u)$ and $\Gamma_k^{ij}(u)$ are smooth functions on $M$. Moreover,
the matrix $F^{ij}(u)$  represents a finite dimensional Poisson structure on
$M$. This  gives a bridge between finite dimensional  and local  Poisson structures.

\bd We say a  local Poisson bracket $\{.,.\}$ in the form  \eqref{loc poiss} admits a dispersionless limit if $F^{ij}(u)=0$ and   $\{.,.\}^{[0]}\neq 0$. In this case  $\{.,.\}^{[0]}$ defines a local
 Poisson bracket on $\lop M$ known as Poisson bracket of
hydrodynamic type. We call it nondegenerate  if $\det \O^{ij}\neq 0$ on an open dense subset of $M$.
\ed

The following theorem, due to Dubrovin and Novikov, relates contravariant metrics on a manifold $M$ to theory of local Poisson brackets on $\lop M$.

\bt\cite{DN}\label{DN thm}
In the notations of formulas \eqref{loc poiss}, if $\{.,.\}^{[0]}$ is  a nondegenerate Poisson brackets of
hydrodynamic type, then the matrix $\O^{ij}(u)$ defines a contravariant flat metric on $M$ and
$\Gamma_{k}^{ij}(u)$ are its  contravariant Christoffel symbols.
\et


We recall the notion of  Dirac reduction of a local Poisson bracket to loop spaces of certain  sub-manifolds. Let us fix a submanifold  $M'\subset M$  of dimension $r$.  We assume $M'$ is defined by the equations $u^\alpha=0$ for $\alpha=r+1,\ldots,n$. We introduce three types of indices; capital letters $I,J,K,\ldots=1,\ldots,n$,
small letters $i,j,k,\ldots=1,\ldots,r$ which parameterize the
submanifold $M'$ and Greek letters
$\alpha,\beta,\gamma,\delta,\ldots=r+1,\ldots,n$.
\bp\cite{mypaper4}\label{dirac fromula} In the notations of equations \eqref{loc poiss}.  Assume  the minor matrix
$F^{\alpha \beta}$ is nondegenerate. Then Dirac reduction is well defined on $\lop {M'}$ and it gives a local Poisson bracket.
If we write the leading terms of the reduced Poisson bracket in the form
\begin{eqnarray}
  \{u^i(x),u^j(y)\}^{[-1]}_{M'} &=& \widetilde{F}^{ij}(u)\delta(x-y), \\
  \{u^i(x),u^j(y)\}^{[0]}_{M'} &=& \widetilde{\O}^{ij} (u)\delta' (x-y)+ \widetilde{\Gamma}_k^{ij} u_x^k \delta (x-y),\\\nonumber
  \{u^i(x),u^j(x)\}^{[k]}_{M'}&=& \widetilde{S}^{ij}_k (u)\delta^{k+1}(x-y)+\ldots,~ k>0.
\end{eqnarray}
Then
\begin{eqnarray}\label{finite dirac red}
\widetilde{F}^{ij}&=& F^{ij}-F^{i\beta} F_{\beta\alpha} F^{\alpha j}
,\\\nn
\widetilde{\O}^{ij}&=& \O^{ij}-\O^{i\beta} F_{\beta\alpha}F^{\alpha
j}+F^{i \beta}F_{\beta \alpha} \O^{\alpha \varphi} F_{\varphi
\gamma} F^{\gamma j}-F^{i\beta} F_{\beta \alpha} \O^{\alpha j},\\\nn
\widetilde{\Gamma}^{ij}_k u_x^k& = & \big(\Gamma^{ij}_k - \Gamma^{i\beta}_k  F_{\beta \alpha} F^{\alpha j} + F^{i \lambda} F_{\lambda \alpha} \Gamma^{\alpha \beta}_k F_{\beta \varphi} F^{\varphi j}-F^{i\beta} F_{\beta \alpha} \Gamma^{\alpha j}_k \big) u_x^k\\\nn & &-\big(\O^{i\beta} - F^{i\lambda} F_{\lambda \alpha} \O^{\alpha \beta} \big)\partial_x(F_{\beta \varphi} F^{\varphi j}),
\end{eqnarray}
while  other higher terms  could be found by solving certain recursive equations.
\ep

\bc\label{cor dirac}
$\widetilde{F}^{ij}$  is the Dirac reduction of the finite dimensional Poisson structure $F^{IJ}$ on $M$ to $M'$. If the entries  $F^{i \alpha}=0$ on $M'$, then the reduced Poisson bracket on $\lop {M'}$  have the same leading terms, i.e.,
 \be  \widetilde{F}^{ij}= F^{ij},~~ \widetilde{\O}^{ij}=\O^{ij},~~
 \widetilde{\Gamma}^{ij}_k = \Gamma^{ij}_k,~~\mathrm{and }~~ \widetilde S^{ij}_k=S^{ij}_k. \ee
\ec

\subsection{From  bihamiltonian structures to Frobenius manifolds}

We use the notations given in section \ref{first-sec} to bring a relations between local bihamiltonian structures   and Frobenius manifolds.

\bd \cite{DFP}
Let $\O_1^{ij}$ and $\O_2^{ij}$ be  two flat contravariant metrics on $M$ with Christoffel symbols $\Gamma_{2k}^{ij}$ and $ \Gamma_{1k}^{ij}$, respectively. Then they form a flat pencil of metrics if  $\O_\lambda^{ij}:=\O_2^{ij}+\lambda \O_1^{ij}$ defines a flat metric on $T^*M$ for generic $\lambda$ and the Christoffel symbols of $\O_\lambda^{ij}$ satisfy  $\Gamma_{\lambda k}^{ij}=\Gamma_{2k}^{ij}+\lambda \Gamma_{1k}^{ij}$. Such   flat pencil of metrics  is
called  quasihomogenous of  degree $d$ if there exists a
function $\tau$ on $M$ such that the vector fields
\begin{eqnarray} \label{tau flat pencil} E&:=& \nabla_2 \tau, ~~E^i
=\O_2^{ij}\partial_{u^j}\tau
\\\nonumber  e&:=&\nabla_1 \tau, ~~e^i
= \O_1^{ij}\partial_{u^j}\tau  \end{eqnarray} satisfy the following
properties
\be  [e,E]=e,~~ \Lie_E \O_2^{ij} =(d-1) \O_2^{ij},~~ \Lie_e \O_2^{ij} =
\O_1^{ij}~~\mathrm{and}~~ \Lie_e\O_1^{ij}
=0.
\ee
 In addition, the  quasihomogenous flat pencil of metrics is called \textbf{regular} if  the
(1,1)-tensor
\begin{equation}\label{regcond}
  R_i^j = {d-1\over 2}\delta_i^j + {\nabla_1}_i
E^j
\end{equation}
is  nondegenerate on $M$.
\ed

The connection between the theory of Frobenius manifolds and flat pencil of metrics is encoded in the following theorem due to Dubrovin.

\bt\cite{DFP}\label{dub flat pencil}
A contravariant quasihomogenous regular  flat pencil of metrics of degree $d$ on a manifold $M$ defines a Frobenius  structure on $M$ of charge $d$.
\et

It is well known that from a Frobenius manifold we always have a flat pencil of metrics but it does not necessarily satisfying the regularity condition \eqref{regcond} \cite{DFP}. Locally, in the coordinates defining equations  \eqref{frob} and \eqref{quasihomog}, the flat pencil of metrics is
found by setting  \begin{eqnarray}\label{frob eqs} \O_1^{ij}&=&\eta^{ij}, \\
\nonumber \O_2^{ij}&=&(d-1+d_i+d_j)\eta^{i\alpha}\eta^{j\beta}
\partial_{t^\alpha}
\partial_{t^\beta} \mathbb{F}.
\end{eqnarray}
This flat pencil of metric is quasihomogeneous  of degree $d$ with $\tau =t^1$. Furthermore, we have
\begin{equation}
E=\sum_i d_i t^i {\partial_{t^i}};~~~e={\partial_{t^{r}}}.
\end{equation}

There is a source of flat pencil of metric within the theory of local bihamiltonian structures.

\bd Two local Poisson brackets  $\{.,.\}_1$ and  $\{.,.\}_2$ on $\lop M$ form  a bihamiltonian structure or they are compatible if the Poisson pencil $\{.,.\}_\lambda:= \{.,.\}_2+\lambda \{.,.\}_1$ is a Poisson bracket  for generic constant $\lambda$. Compatible Poisson brackets $\{.,.\}_1$ and  $\{.,.\}_2$ form an an exact Poisson pencil if there exists a vector field $X$ such that
\be \{.,.\}_1 = \Lie_X \{.,.\}_2;~~\Lie_X \{.,.\}_1= 0.
\ee
In this case we call $X$  Liouville vector field.
\ed

 For recent developments about the theory of exact Poisson pencil see \cite{FalLor} and \cite{LPM}.

    Let us fix  compatible local  Poisson brackets $\{.,.\}_2$ and  $\{.,.\}_1$  on $\lop M$ and write their  leading  terms   in the form
\begin{eqnarray}\label{bihamileading}
  \{u^i(x),u^j(y)\}^{[-1]}_\alpha &=& F^{ij}_\alpha(u(x))\delta(x-y), \alpha=1,2\\\nn
  \{u^i(x),u^j(y)\}^{[0]}_\alpha &=& \O^{ij}_{\alpha}(u(x)) \delta' (x-y)+ \Gamma_{\alpha k}^{ij}(u(x)) u_x^k \delta (x-y).
\end{eqnarray}
 Suppose that $\{.,.\}_1$ and $\{.,.\}_2$ admit a dispersionless limit (we also say the bihamiltonian structure admits a dispersionless limit). In addition, assume the corresponding Poisson brackets of hydrodynamics type are nondegenerate as well as the dispersionless limit of $\{.,.\}_\lambda$ for generic $\lambda$. Then, using theorem \ref{DN thm}, the matrices $\O_1^{ij}$ and $\O_2^{ij}$ define a  \textbf{flat pencil of metrics} on $M$.

\section{Nilpotent elements of semisimple type }

In this section, we collect properties of the so called  distinguished nilpotent elements of semisimple type in simple Lie algebras. Then  we  derive important identities needed to prove our main results.
\subsection{Background}
We fix a complex simple Lie algebra $\g$ of rank $r$. We refer to the type of $\g$ by $Z_r$.  For $g\in \g$, let  $\mathcal O_g$ denotes  the orbit of $g$ under the adjoint group action.  The element  $g$ is called nilpotent if $\ad_g$ is  nilpotent in $\textrm{End}(\g)$ and it is called regular if $\dim \g^{g}=r$. Any simple Lie algebra contains regular nilpotent elements.

We fix a nilpotent element  $L_1$  in $\g$ (later, we will assume  it is  distinguished). Let  $A:=\{L_1,h,f\}\subseteq \g$ be an associated  $sl_2$-triple satisfying the relations \eqref{sl2:relation1}. It follows from representation theory of $sl_2$-algebra that  the eigenvalues of $\ad_h$ are integers and half integers. Consider Dynkin grading associated to $L_1$
\begin{equation}\label{grad1}
\g=\bigoplus_{i\in {1\over 2}\mathbb{Z}} \g_i;~~~~\g_i:=\{g\in\g: \ad_h g= i g\}.
\end{equation}

 We retrieve from \cite{COLMC} the following definitions concerning nilpotent orbits and their classification.  If $L_1$ is  regular, then   $\mathcal{O}_{L_1}$ is called  regular nilpotent orbit, and it is equal to the set of all regular nilpotent elements in  $\g$. The  nilpotent orbit  $\mathcal{O}_{L_1}$    is called {\bf distinguished}, and hence also $L_1$, if  $\mathcal{O}_{L_1}$ has no representative in a proper Levi subalgebra of $\g$. It turns out that $L_1$ is distinguished iff $\dim \g_0=\dim \g_1$. Moreover, if $L_1$ is distinguished, then the eigenvalues of $\ad_h$ are all integers.  The regular nilpotent orbit in $\g$ is  distinguished.

 Distinguished nilpotent orbits, along with other nilpotent orbits, are classified by using weighted Dynking diagrams.  In the case $\g$ is an exceptional  Lie algebra,  distinguished nilpotent orbits are listed in the form $Z_r(a_i)$ where  $i$ is the number of vertices of weight 0 in the corresponding weighted  Dynkin diagram.  If there is another orbit of the same number $i$ of 0's, then the notation $Z_r(b_i)$ is used.  For all simple Lie algebras, the type of the regular nilpotent orbit is  $Z_r(a_0)$.

 In case $\g$ is a classical  Lie algebra, nilpotent orbits  are also classified through partitions of the dimension of the fundamental representation of $\g$. In this article, by $B_{2m}(a_m)$, we refer to the distinguished nilpotent orbit  corresponding to the partition $[2m+1,2m-1,1]$ when the Lie algebra  $\g$ is  $so_{4m+1}$ (type $B_{2m}$). While, as usual in the literature,   $D_{2m}(a_{m-1})$  denotes  the distinguished nilpotent orbit corresponding to the partition $[2m+1,2m-1]$ when   $\g$ is  $so_{4m}$ (type $D_{2m}$).

From \cite{Elash}, we recall the following definition and properties. The nilpotent element  $L_1$ is of  semisimple type, and so its orbit, if there exists an element $g$ of the minimal eigenvalue of $\ad_{h}$ such that $L_1+g$ is semisimple. In this case $L_1+g$ is called  a cyclic element. If  $L_1$ is also distinguished then $L_1+g$ will be regular.   The list of distinguished nilpotent elements of semisimple types  are (\textit{idid}, Lemma 5.1 and see the appendix of \cite{DelFeher}):
\begin{enumerate}
    \item All regular nilpotent orbits  in simple Lie algebras (those of type $Z_r(a_0)$)
    \item Subregular nilpotent orbits $F_4(a_1)$, $E_6(a_1)$, $E_7(a_1)$ and $E_8(a_1)$.
    \item Nilpotent orbits of type $B_{2m}(a_m)$ and  $D_{2m}(a_{m-1})$.
    \item Nilpotent orbits of type $F_4(a_2)$,  $F_4(a_3)$, $E_6(a_3)$,   $E_7(a_5)$, $E_8(a_2)$,$E_8(a_4)$, $E_8(a_6)$ and $E_8(a_7)$.
\end{enumerate}

 From now on,    we assume that $L_1$ is  a distinguished nilpotent element of semisimple  type  and we refer to its type by $Z_r(a_s)$.    Let    $\eta_r$ denotes the maximal eigenvalue of $\ad_h$. Thus, we can (and will) fix an element $K_1\in \g_{-\eta_r}$ such that the cyclic element  $\Y_1:=L_1 +K_1$ is regular semisimple.

 In what follows we give  a general setup associated to the  cyclic element $\Y_1$. It was initiated by  Kostant for the case of  regular nilpotent elements \cite{kostBetti} and obtained for  distinguished nilpotent elements of semisimple type in  \cite{DelFeher}.    Let $\h':=\g^{\Y_1}$ be the Cartan subalgebra  containing  $\Y_1$ which is known as the opposite Cartan subalgebra. Then  the adjoint group element $w$ defined by
\begin{equation} \label{qcoxeter}
 w:=\exp {2\pi \mathbf{i}\over {\eta_r+1}}\ad_h
\end{equation}
acts on $\h'$ as a representative of a regular cuspidal conjugacy class $[w]$  in the Weyl group $\W(\g)$ of $\g$ of order ${\eta_r+1}$. We recall that   a conjugacy class $[w']\subset \W(\g) $  is called
cuspidal (resp. primitive) if $\det(w'-I)\neq 0$ (resp. $\det(w'-I)=\det \mathbb K$, $\mathbb K$ is the Cartan matrix of $\W(\g)$).  Also, $[w']$ is called
regular  if $w'$ has an
eigenvector not fixed by any non-identity element in $\W(g)$ (see \cite{springer} for the  classification of regular conjugacy classes). We emphasize that the  results in this article depend on the nilpotent orbit  $\mathcal{O}_{L_1}$ and not on the particular representative $L_1$ of  $\mathcal{O}_{L_1}$.

\subsection{Normalization and identities}
The  element $\Y_1$ is an eigenvector of $w$ of eigenvalue $\epsilon=\exp \frac{2\pi \mathbf{i}}{\eta_r+1}$. We define the multiset  $\exx(L_1)$ which  consists of natural numbers  $\eta_i$, $i=1,\ldots,r$ such that  $\epsilon^{\eta_i}$'s is an eigenvalue of the action of  $w$ on $\h'$.   We call  $\exx(L_1)$  the exponents of the nilpotent element $L_1$. When  $L_1$ is a regular nilpotent element, $\exx(L_1)$ equals the  exponents $\exx(\g)$ of the Lie algebra  \cite{kostBetti}. In table 1, we list  elements of  $\exx(L_1)$  in the second column. We calculated them  by combining the results of  \cite{DelFeher},\cite{Elash} and \cite{springer}. Note that $\exx(\g)$ is listed in table \ref{we} as $\exx(L_1)$ when $L_1$  is of  type $Z_r(a_0)$. We  denote throughout this article,  the elements of $\exx(L_1)$ by $\eta_i$ and elements of $\exx(\g)$ by $\nu_i$ and  we assume they are given in a non-decreasing order, i.e.,
\be \eta_1\leq \eta_2\leq \ldots \leq \eta_r ~\textrm{and} ~\nu_1\leq \nu_2\leq \ldots \leq \nu_r.\ee

 The following lemma summarize an important relation between $\exx(\g)$ and $\exx(L_1)$.

\bl \label{relat-expon} For $i=1,\ldots,r$,
\be \eta_{i}+\eta_{r-i+1}=\eta_r+1.
\ee
Moreover, there exists a unique non-negative integer $\mu_i$ such that $\nu_i-\mu_i(\eta_r+1)$ belongs to $\exx(L_1)$. Furthermore, the multiset formed by the numbers
 $\nu_i-\mu_i(\eta_r+1)$  equals  the multiset $\exx(L_1)$. In addition, the number of $\mu_i's$ which are zero equals  $r-s$.
 \el
 \begin{proof}
The proof is obtained  by examining the multisets $\exx (L_1)$ and $\exx(\g)$  for each nilpotent orbit listed in table \ref{we}.
 \end{proof}

\bx  In case
 $L_1$ is of type $E_7(a_5)$. Then  $\eta_7=5$ and  the values of $u_i's$ are given in the following table

\begin{center}
    \begin{tabular}{|c|c|c|c|c|c|c|c|}
      \hline
      $i$ & 1& 2& 3& 4& 5& 6& 7 \\\hline
      $\nu_i$ & 1 & 5 & 7 & 9 & 11 & 13& 17 \\\hline
     $\mu_i$ & 0 & 0 & 1 & 1 & 1 & 2 & 2\\\hline
      $\nu_i-\mu_i(\eta_7+1)$& 1 & 5 & 1 & 3 & 5 & 1 & 5  \\
      \hline
    \end{tabular}
    \end{center}
The last row is just the elements of $\exx(L_1)$ (not in order).
\ex

We  keep the notations $\mu_i$, $i=1,\ldots,r$ for the non-negative numbers introduced in the last lemma.   Many formulas below  depend on these numbers.
We list them in the fourth column of table \ref{we} using conventional notation for repetitions. For example, we write $[0^2,1^3,2^2]$ instead of $[0,0,1,1,1,2,2]$.

Let  $\Y_1,\Y_2,\ldots, \Y_{r}$ be a basis of  $\h'$ of eigenvectors of $w$ such that  $w(\Y_i)=\epsilon^{\eta_i} \Y_i$. Then $\Y_i$ has the form
\be \Y_i=L_{i}+K_{i}; ~~L_{i} \in \g_{\eta_i},~ K_{i}\in \g_{\eta_i-(\eta_r+1)},~~L_i\neq 0\neq K_i,~~i=1,\ldots,r.\ee

\begin{table}
\begin{center}
\small
\begin{tabular}{|c|c|c|c|}\hline
 &  \multicolumn{2}{|c|}{$\ew(L_1)$}&  \\
\hline
$Z_r(a_s)$ &  $\eta_1\leq\ldots\leq \eta_{r}$&$\eta_{r+1}\leq \ldots\leq \eta_{n}$ &$[\mu_1,\ldots, \mu_r]$ \\ \hline

$A_{r}(a_{0})$ &  $1,2,\ldots,r$& - & $[0^{r}]$\\ \hline
$B_{r}(a_{0})$ &  $1,3,\ldots,2r-1$& - & $[0^{r}]$\\\hline
$B_{2m}(a_{m})$ &  $1,1,3,3,\ldots,2m-1,2m-1$&$1,2,\cdots,m-1;m-1,$ & $[0^{m+1},1^{m-1}]$\\
& &  $m;m,m+1,\ldots,2m-2$&\\
\hline
$C_{r}(a_{0})$ &  $1,3,\ldots,2r-1$& - & $[0^{r}]$\\\hline
$D_{r}(a_{0})$ &  $1,3,\ldots,r-1;r-1,r-3,\ldots,2r-3$& - & $[0^{r}]$\\\hline
$D_{2m}(a_{m-1})$ &  $1,1,3,3,\ldots,2m-1,2m-1$&$1,2,\cdots,2m-2$& $[0^{m+1},1^{m-1}]$\\\hline

$E_6(a_0)$& $1,4,5,7,8,11$  &-& $[0^{n}]$\\
$E_6(a_1)$& $1,2,4,5,7,8$ & $3,5$& $[0^{n-1},1]$\\
$E_6(a_3)$& $1,1,2,4,5,5$  &$1,2,2,3,3,4$& $[0^3,1^3]$\\\hline
$E_7(a_0)$ & $1,5,7,9,11,13,17$  & -& $[0^7]$\\
$E_7(a_1)$ & $1,3,5,7,9,11,13$  &  $5,8$& $[0^6,1]$\\
$E_7(a_5)$ & $1,1,1,3,5,5,5$  & $1,1,1,2,2,2,2,$& $[0^2,1^3,2^2]$\\
 &   & $ 3,3,3,3,4,4,4$& \\ \hline
 $E_8(a_0)$ & $1,7,11,13,17,19,23, 29$  & - &$[0^8]$\\

$E_8(a_1) $ & $1,5,7,11,13,17,19,23$  & $9,14$&$ [0^7,1]$\\
$E_8(a_2) $ & $1,3,7,9,11,13,17,19$  & $5,8,11,14$& $[0^6,1^2]$\\
$E_8(a_4) $ & $1,2,4,7,8,11,13,14$  &$3,5,5,7,7,9,9,11$  &$[0^5,1^3]$\\
$E_8(a_5) $ & $1,1,5,5,7,7,11,11$   &$1,2,3,4,5,5,6,6,7,8,9,10$ &$[0^3,1^4,2]$\\
$E_8(a_6) $ & $1,1,3,3,7,7,9,9$  &$1,2,3,3,3,4,4,4, 5,5,5,6,6,6,7,8$ &$[0^2,1^4,2^2]$\\

$E_8(a_7) $ & $1,1,1,1,5,5,5,5$ &$1, 1, 1, 1, 1, 1,2, 2, 2, 2, 2, 2, 2, 2, 2, 2 $ &$[0,1^2,2^2,3^2,4]$\\
&  & $3,3,3,3,3,3,3,3,3,3,4,4,4,4,4,4$&\\

\hline
$F_4(a_0)$ & 1,5,7,11 & -&$[0^4]$\\
$F_4(a_1)$ & 1,3,5,7 & 2,5&$[0^3, 1]$\\
$F_4(a_2)$ & 1,1,5,5  &1,2,3,4& $[0^2, 1^2]$\\
$F_4(a_3)$ &  1,1,3,3 &1,1,1,1,2,2,2,2& $[0,1^2,2]$\\\hline
$G_{2}(a_{0})$ &  $1,5$& - & $[0^{2}]$\\ \hline
& \multicolumn{1}{|c|}{$\exx(L_1)$}& $\overline\exx(L_1)$&\\\hline
\end{tabular}
\caption{Exponents and  weights of distinguished nilpotent elements of semisimple type}
\label{we}
\end{center}
\end{table}

We normalized  the invariant nondegenerate bilinear form $\bil . .$ on $\g$ such that $\bil{L_1}{f}=1$. Then the following lemma is valid.
\bl\label{opp Cartan prop}
The matrix $T_{ij}:=\bil {\Y_i}{\Y_j}$ is nondegenerate and  antidiagonal with respect to $\exx(L_1)$, i.e.,  \[T_{ij}=0, ~{\rm if}~\eta_i+\eta_j\neq \eta_r+1.\]
Moreover, The elements  $\Y_i,~i>1$  can be normalized such that
\be
\bil{\Y_i}{\Y_j}=(\eta_r+1)\delta_{i+j,r+1}.
\ee
\el
\begin{proof}
The first part follows from the fact that the restriction of $\bil . .$ to a Cartan subalgebra is nondegenerate. Therefore, for any element $\Y_i$ there exists an element  $\Y_j$ such that $\bil {\Y_i}{\Y_j}\neq 0$. But for  the Weyl group element  $w$ defined in \eqref{qcoxeter}, we have the equality
\[\bil {\Y_i}{\Y_j}=\bil {w \Y_i}{w \Y_j}=\exp {2(\eta_i+\eta_j)\pi \mathbf{i} \over {\eta_r+1}}\bil {\Y_i}{\Y_j} \]
which forces  $\eta_i+\eta_j=\eta_r+1$ in case $\bil {\Y_i}{\Y_j}\neq 0$. For the
second part of the lemma, recursively, we  can define  a change of basis with linear combination   upon the elements $\Y_i$ which have  the same eigenvalue such that the matrix $T_{ij}$ transform to the anti-diagonal form  $(\eta_r+1)\delta_{i+j,r+1}$.

\end{proof}

 We  assume from now on that  the basis $\Y_i$ of $\h'$ are normalized and satisfy the hypothesis of the previous lemma. Then we get the following identities.

\bc \label{silver}
\be
\bil{L_i}{K_{j}}=\eta_j\delta_{i+j,r+1},~~i,j=1,\ldots,r.
\ee
\ec
\begin{proof}
Recall that
 \be \Y_i=L_{i}+K_{i}; ~~L_{i} \in \g_{\eta_i},~ K_{i}\in \g_{\eta_i-(\eta_r+1)}.\ee
Using the relation $0=[\Y_i,\Y_j]=[L_i,K_j]+[K_i,L_j]$ with the invariant bilinear form yields
\be 0=\bil h {[L_i,K_j]+[K_i,L_j]}= (\eta_i )\bil {L_i}{K_j}+ (\eta_i-(\eta_r+1)) \bil {K_i} {L_j}. \ee
This equation  with the normalization   $\bil {\Y_i}{\Y_j}=\bil {L_i}{K_j}+\bil {K_i}{L_j}=(\eta_r+1)\delta_{i+j,r+1}$ lead to the required identity.
\end{proof}
\bc\label{Gold}
\[\bil {[K_1,L_j]}{\ad_f L_i}=\eta_i \eta_j\delta_{i+j,r+1}, ~~i,j=1,\ldots,r.\]

\ec
\begin{proof}
The identity $[\Y_1,\Y_j]=0$ leads to $[L_1,K_j]=-[K_1,L_j]$.
Then
\begin{eqnarray}
\bil {[K_1,L_j]}{\ad_f L_i}&=&-\bil {[L_1,K_j]}{[f,L_i]}=\bil {K_j}{[L_1,[f,L_i]]}=\bil {K_j}{[L_i,[f,L_1]]}\\\nn
&=& -\bil {K_j}{[L_i,h]}=\eta_i \bil {K_j}{L_i}=\eta_i \eta_j\delta_{i+j,r+1}.
\end{eqnarray}
\end{proof}

The commutators $[\Y_i,\Y_j]=0$ imply  that the set  $\{L_{1},\ldots,L_{r}\}$ generates a commutative subalgebra of $\g^{L_1}$. We consider the restriction of the  adjoint representation to the $sl_2$-subalgebra $\A$ generated by $\{L_1,h,f\}$. Then  the vectors  $L_{i}$ are  maximal weight vectors of  irreducible $\A$-submodules $\V_i$  of dimension  $2\eta_i+1$. We set $n=\dim \g^{L_1}$ and we fix the following  decomposition of $\g$ into irreducible $\A$-submodules
\begin{equation}\label{decompo}
\g=\bigoplus_{j=1}^{n} \V_j,~~\dim \V_j=2\eta_j+1,~ L_j\in \V_j, \ad_{L_1} L_j=0,~ \ad_h L_j=\eta_j L_j.
\end{equation}
  Note that, for convenience, we extend the notation $L_j$ to  cover all  maximal eigenvectors, i.e., $L_j$'s form a basis for $\g^{L_1}$.  The numbers  $\eta_1,\ldots,\eta_{n}$ are given in  table \ref{we} as the collection of the numbers in the second and fourth columns.  We refer to them as the weights of the nilpotent element $L_1$. We could not find them in the literature and we had to calculate them explicitly. See \cite{mypaperw} for a procedure to find the weights of a distinguished nilpotent element and the calculation for the nilpotent element of type $D_{2m}(a_{m-1})$.
 After calculating the weights, we observe the following.

 \bc
  $n=r+2 \sum \mu_i$.
 \ec
 Let $\overline \exx(L_1)$ denotes the multiset consisting of the numbers $\eta_i$, $i=r+1,\ldots,n$ and assume they are given in non-decreasing order, i.e.,
\[\eta_{r+1}\leq \eta_{r+1}\leq \ldots \leq \eta_{n}.\]
Then from table \ref{we}, we get

\bc \label{Ecompl} $\eta_{r+i}+\eta_{n-i+1}=\eta_r$ for $i=1,\ldots, n-r$.
\ec

We use the fact that $\g^f$ is the dual of $\g^{L_1}$ under  $\bil . .$ \cite{wang}  to fix a basis $\ga_i$ of  $\g^f$ such that
\be
\bil {\ga_i} {L_j}=\delta_{ij},~ i=1,\ldots,n.
\ee
Then   $\ad_h\ga_i=-\eta_i \ga_i$. Let us  introduce the following  basis for  $\bigoplus_{i\leq 0} \g_i$
\be
\ga_i, \ad_{L_1} \ga_i,\ldots, {1\over \eta_i!} \ad_{L_1}^{\eta_i} \ga_i,~~ i:=1,\ldots,n,
\ee
and similarly a basis for $\bigoplus_{i\geq 0} \g_i$
\be
L_i, \ad_{f} L_i,\ldots, \ad_{f}^{\eta_i} L_i,~~ i:=1,\ldots,n.
\ee
\bl
 \be
 \bil {{1\over I!} \ad_{L_1}^{I} \ga_i}{ad_f^J L_j}=(-1)^I \bin {\eta_i}{I} \delta_{ij}\delta^{IJ};~I=0,1,\ldots,\eta_i; J=0,1,\ldots ,\eta_j.
 \ee
\el
\begin{proof}
For $I=J=1$, we get
\be \bil {\ad_{L_1}  \ga_i}{ad_f L_j}=-\bil {\ga_i}{ad_{L_1}  \ad_f L_j}=\bil {\ga_i}{[L_j,h]}=-\eta_i \delta_{ij}.
\ee
Hence, by induction for $I>1$,
\begin{eqnarray}
\bil {{1\over I!} \ad_{L_1}^I \ga_i}{ad_f^I L_j}&=&\bil {{1\over I!} \ad_{L_1}^{I-1} \ga_i}{[ad_f^{I-1} L_j,h]}\\\nonumber
&=& -{{\eta_j-I+1}\over I}\bil {{1\over (I-1)!} \ad_{L_1}^{I-1} \ga_i}{ad_f^{I-1} L_j}=(-1)^I \bin {\eta_i}{I} \delta_{ij}.
\end{eqnarray}
Suppose $I>J$. Then we can recursively equate  the value $ \bil { \ad_{L_1}^{I} \ga_i}{ad_f^J L_j}$ to constant multiplication of the zero valued $ \bil { \ad_{L_1}^{I-J-1} \ga_i}{ad_f L_j}$.
\end{proof}

\bc
$\gamma_r=K_1$.
\ec
\begin{proof}
Recall that $K_1\in \g_{-\eta_r}$. It follows from the Dynkin grading that $K_1\in \g^f$. Then,  for $j\leq r$, it  follows from corollary \ref{silver} that $\bil {K_1}{L_j}=\delta_{jr}$. While for $j>r$, we get from Dynkin grading and the fact that $\eta_j<\eta_r$ that $\bil {K_1}{L_j}=0$. Thus by construction $\gamma_r=K_1$.
\end{proof}

\section{The space of common equilibrium points}

In this section, we  fix  Slodowy slice $Q$ as a transverse subspace to the orbit space of $L_1$. We  discuss the integrability of the transverse Poisson structure at $L_1$ of  Lie-Poisson structure on $\g$ which leads to the definition of the space of common equilibrium points $N$. Then, we  will introduce special coordinates on $Q$ and give alternative  definitions  for $N$.

\subsection{Background} \label{background1}
Let us define   the gradient $\nabla H: \g\to \g$ for a function $H$ on $\g $ by
\be \label{gradient}
\frac{d}{dt}_{|_{t=0}}H(g+tv)=\bil {\nabla H(g)}{v},~\forall g, v\in \g.\ee
We fix the following standard  compatible Poisson structures on $\g$ which consists of the frozen Lie-Poisson structure $B_1^\g$ and the standard Lie-Poisson structure $B_2^\g$. We denote their Poisson brackets by   $\{.,.\}_1^\g$ and $\{.,.\}_2^\g$, respectively. For any two functions $H$ and $G$ on $\g$, and  $v\in T_g^* \g\cong\g$, we set
\begin{eqnarray}
 \label{bih:stru on g}
\{H,G\}_1^\g(g)&=&\bil {[ \nabla G(g),\nabla H(g
)]}{K_1}; ~~B_1^\g(v)= [K_1,v],\\\nonumber
\{H,G\}_2^\g(g)&=&\bil{[ \nabla G(g),\nabla H(g
)]} {g}; ~~B_2^\g(v)=[ g,v].
\end{eqnarray}
 We   use $B_i^\g$, $i=1,2$ to refer  to both  the Poisson structures (tensors) and the corresponding Poisson brackets. Then the  Hamiltonian vector field $\chi_H$ of  a function $H$ under $B_2^\g$  at a point $g\in \g$  is defined by
\be \chi_H(g):=-\ad_{\nabla H(g)} g=[ g,{\nabla H(g)}].\ee
It is known that the symplectic leaf through $g\in \g$  coincides with the adjoint orbit  $\mathcal O_g$  and invariant polynomials under the adjoint group action are global Casimirs of $B_2^\g$.

 Using  Chevalley's theorem, we fix  a complete  system of homogeneous generators $P_{1},\ldots,P_{r}$ of the  ring of invariant polynomials under the adjoint group action.  We  assume that degree $P_i$ equals $\nu_i+1$. These generators  give a complete set of  global Casimir functions of  $B_2^\g$. In particular,
\be \label{center grad} \nabla P_i(g)\in \g^g,~ \forall g\in \g,~i=1,\dots,r.\ee
 Moreover, the functions $P_i(g+\lambda K_1)$ form a complete set of independent global Casimirs of the Poisson pencil $B_\lambda^\g:=B_2^\g+\lambda B_1^\g$ for any  $\lambda\in \mathbb C$ \cite{bolv1}.

Define Slodowy slice  $Q$ to be the affine space
\be Q:=L_1+\g^f. \ee
Then  $Q$  is a transverse subspace to the symplectic leaf $\mathcal O_{L_1}$ of $B_2^\g$ through  $L_1$.  The following proposition is a special version of theorem \ref{psred} stated below.

\bp \cite{mypaper4}
The space $Q$ inherits compatible Poisson structures   $B_1^Q$, $B_2^Q$ from $B_1^\g$, $B_2^\g$, respectively. Moreover, $B_2^Q$ is the transverse Poisson structure at $L_1$ of  Lie-Poisson structure  $B_2^\g$. Furthermore,  for any $\lambda\in \mathbb C$,  $B_\lambda^Q:=B_2^Q+\lambda B_1^Q$ can be obtained from $B_\lambda^\g$   using  Dirac reduction.
\ep

Let $\overline P_i^0$ denotes the restriction of the invariant polynomial $P_i$ to $Q$.
Since $B_\lambda^Q$ can be obtained by Dirac reduction, we have the following standard consequence.
\bp
For $\lambda\in \mathbb C$, $\overline P_1^0(q+\lambda K_1),\ldots, \overline P_r^0(q+\lambda K_1)$  form a complete set of independent  Casimirs of the Poisson pencil $B_\lambda^Q$.
\ep

Following  the argument shift method (\cite{bolv}, \cite{Mfomenko}), we consider  the family of functions
\be
\mathbf{F}:=\cup_{\lambda\in \overline \mathbb{C}}\{P'_\lambda: P'_\lambda ~\mathrm{is~ a~ Casimir~ of} ~B_\lambda^Q\}.
\ee
This family commutes pairwise with respect to both Poisson brackets (\cite{bolv}, section 1.3).
Let us
consider the coefficient $\overline P_i^j$ of Taylor expansions
 \be \label{arg:sht}
 \overline P_i^0(q+\lambda K_1)=\sum_{j\geq 0} \lambda^{j} \overline  P_i^j(q), ~q\in Q.
 \ee
 Then the functions $\overline P_i^j$ functionally generate $\mathbf F$. Moreover, $\overline P_i^0$ are Casimirs of $B_2^Q$, the highest non-constant term  $\overline P_i^{\varrho_i}$ are Casimirs of $B_1^Q$, and  all functions $\overline P_i^j$ are in involution with respect to both Poisson structures. In proposition \ref{norm-coord}, we will  show that $\varrho_i=
\mu_i$.

The main propose for applying argument shift method is to show that $\mathbf F$ contains enough number  of functionally independent  functions in order to get a completely integrable system for $B_2^Q$. We explored this problem in  \cite{mypaper5}  for arbitrary nilpotent elements in $\g$ and we proved the following theorem

\bt \label{mainthm2} \cite{mypaper5}
Suppose  $L_1$ belongs to one of  the following  distinguished nilpotent orbits  of  semisimple type: $D_{2m}(a_{m-1})$, $B_{2m}(a_m)$,  $F_4(a_2)$,  $E_6(a_3)$,   $E_8(a_2)$ and $E_8(a_4)$.
Then the set of all   functions ${\overline P_i^j}$ result from the expansion \eqref{arg:sht} are  functionally independent and   form a polynomial  completely integrable system under $B_2^Q$.
\et

In what follows a point $q\in Q$ is generic if $\rk\, B_2^Q(q)=n-r$. From \cite{bolv1} we get the following theorem

\bt \cite{bolv1}  \label{bolsinov1}
The family $\mathbf F$ is complete (contains a completely integrable system) if and only if, at a generic point $q\in Q$, $\rk\, B_{\lambda}^Q(q)=\rk\, B_{\zeta}^Q(q)$    for all  $\lambda, \zeta\in \overline \mathbb C$.
\et

We are concern about the space of common equilibrium points $N$ of the family $\mathbf F$ which is defined by
 \be
 N:=\{q\in Q: B_\lambda^Q(dP')(q)=0,~ \forall P'\in \mathbf F, \lambda \in \overline{\mathbb C}\}.
 \ee
The following theorem gives an equivalent definition.
\bt\cite{bolv1} \label{bolsinov2}
A point $q\in Q$ is a common equilibrium point  if and only if $\ker B_{\lambda}^Q(q)=\ker B_{\zeta}^Q(q)$    for all generic  $\lambda, \zeta\in \overline \mathbb C$.
\et

Equivalently, for $q$ to be in $N$, it is sufficient to require that  the kernel of just two generic brackets at $q$ coincides, i.e., $\ker B_\lambda^Q(q)=\ker B_{\zeta}^Q(q)$ with $\lambda\neq \zeta$ \cite{bolv1}.

\subsection{Special coordinates}

Let us consider  the adjoint quotient map
 \be \Psi:\g\to \mathbb{C}^r,~~~ \Psi(g)=(P_1(g),\dots,P_{r}(g)).\ee
 Kostant  proved in  \cite{kostpoly} that the rank of  $\Psi$ at $g$ equals $r$ if and only if $g$ is a regular element in $\g$ and it is known that the set of regular element is open and dense in $\g$.  Later, Slodowy proved that the rank of $\Psi$ is $r-1$ at subregular nilpotent elements \cite{sldwy2}. Finally,  Richarson   \cite{richard} obtained the ranks of $\Psi$ at distinguished  nilpotent elements  except for the nilpotent elements of type $E_8(a_2)$. Results in  this section are  build on and inspired by the  articles mentioned in this paragraph.

We fix a basis $e_0,e_1,e_2,\ldots$ for $\g$ such that \begin{enumerate}
    \item The  elements $e_0,e_1,\ldots,e_{n+r}$ are  $K_r,L_1,L_2,\ldots,L_n,K_1,K_2,\ldots,K_{r-1}$, respectively. Recall that  $\Y_i= L_i+K_i$ are normalized according to lemma \ref{opp Cartan prop}.
    \item $\bil{e_i}{\Y_1}\neq 0$ if and only if $i=0$ or $i=r$.
 \end{enumerate}
It is not hard to show that such a basis exists. Let us  define on $\g$ the linear coordinates
\be \label{coord on g}
z^i(g)=\bil {e_i}{g}, ~~i=0,1,2,\ldots .
\ee
 Then, by definition, $\nabla H=\sum {\partial H\over \partial z^i}e_i$ for any function $H$ on $\g$. Note that the rank of $\Psi$ at $g$  equals  the dimension of the vector space generated by $\nabla P_i(g)$. In particular, since  $\Y_1$ is regular,  the gradients  $\nabla P_i(\Y_1)$ are linearly independent and form a basis for the opposite Cartan subalgebra  $\h'$.  We use these remarks in the following lemma.

\bl \label{first norm}
  The matrix  with entries  ${\partial P_i \over \partial z^j}(\Y_1)$, $i,j=1,\ldots,r$, is non-degenerate. Moreover, $P_i$ have the following form
\be\label{respoly1} \label{respoly2}
P_i= R_i^1+R_i^2+R_i^3
\ee where
\be \label{respoly21}R_i^1=\sum_{a(\eta_r+1)=\nu_i-\eta_r}\theta_{i,a}(z^r)^{a+1}(z^0)^{\nu_i-a},~~ R_i^2=\sum_{a=0}^{\nu_i-1}\sum_{j=1}^{r-1}  c_{i,j,a}(z^r)^{a}(z^0)^{\nu_i-a}(z^j+z^{j+n}),~
\frac{\partial R_i^3} {\partial z^k}(\Y_1)=0, \forall k.\ee Here,  $
c_{i,j,a}$ and $ \theta_{i,a}$ are complex numbers.
\el

\begin{proof}
Since $\nabla P_i(\Y_1)\in \g^{\Y_1}=\h'$ and  $\h'$ has basis $\Y_i=L_i+K_i$, we get
\be
\nabla P_i(\Y_1)=\sum_{j=1}^{r} C_{i,j} \Y_j=\sum_{j=1}^{r} C_{i,j} ( L_j+ K_j)=C_{i,r}(e_0+e_r)+\sum_{j=1}^{r-1}C_{i,j}(e_j+e_{n+j}).
\ee
Hence
\be
C_{i,j}=\left\{
  \begin{array}{ll}
   \frac {\partial P_i}{\partial z^j}(\Y_1)=\frac{\partial P_i}{ \partial z^{j+n}}(\Y_1)
   ,    & 0<j<r\hbox{;} \\
   \\
     \frac{\partial P_i}{ \partial z^r}(\Y_1)=\frac{\partial P_i}{ \partial z^{0}}(\Y_1), & j=r\hbox{;}
  \end{array}
\right.
\ee
and ${\partial P_i\over \partial z^j}(\Y_1)=0$ for other values of $j$. By  definition of the  coordinates and corollary \ref{silver},   $z^j(\Y_1)$ are all zero except $z^r(\Y_1)=1$ and $z^0(\Y_1)=\eta_r$. For $0<j<r$, imposing the condition  ${\partial P_i\over \partial z^j}(\Y_1)\neq 0$ and using the homogeneity of $P_i$, we find  that   ${\partial P_i\over \partial z^j}$ must contain  the polynomial
\be
 \sum_{a=0}^{\nu_i-1} c_{i,j,a} (z^r)^{a}(z^0)^{\nu_i-a-1}, ~~c_{i,j,a}\in \mathbb{C}.
\ee
This gives the formula for $R_i^2$. Note that ${\partial R_i^2\over \partial z^r}(\Y_1)=0$ since $z^j(\Y_1)=0$ for $j\neq 0$ and $j\neq r$. Thus, for ${\partial P_i\over \partial z^r}(\Y_1)$ to  be nonzero, $P_i$ must contain   terms of the form $\Xi_{i,a}=(z^r)^{a+1} (z^0)^{\nu_i-a}$. But then $a$ is constrained by the identity \be {\partial \Xi_{i,a}\over \partial z^r}(\Y_1)=(a+1) (\eta_r)^{\nu_i-a}={\partial \Xi_{i,a}\over \partial z^{0}}(\Y_1)=({\nu_i-a}) (\eta_r)^{\nu_i-a-1}.
\ee
This leads to  the formula for $R_i^1$. The condition on $R_i^3$ is a direct consequence from our analysis.
Finally,  the non-degeneracy condition  follows from the  fact that the vectors  $\nabla P_i(\Y_1)$ are a basis for $\h'$.
\end{proof}

For  Slodowy slice $Q$, we observe that $z^{0}(q)=\bil {K_r}{L_1}=\eta_r\neq 0$  for every $q\in Q$ and $(z^1,\ldots,z^{n})$ define global coordinates on $Q$.   The value of these coordinates at  $\Y_1\in Q$ are  $z^i=\delta^{ir}$.  We  set degree $z^i$ equals $\eta_i+1$ and recall  the following quasihomogeneity theorem  due to Slodowy.

\bt (\cite{sldwy2}, section 2.5) \label{slodowy}  The restriction $\overline P_i^0$ of   $P_i$  to $Q$ is quasi-homogeneous polynomial of  degree $\nu_i+1$.
\et

This theorem leads to the following refinement of the last lemma.

\bp \label{norm-coord}
 The  restrictions $\overline P_i^0$ of the invariant polynomials $P_i$ to $Q$ in the coordinates $(z^1,\ldots,z^n)$  take the form
 \be  \label{norm-coord eq}
\overline P_i^0(z^1,\ldots,z^n)=\sum_{\nu_i-\eta_j=\mu_i(\eta_r+1)} \widetilde{c}_{i,j}(z^r)^{\mu_i} z^j+\overline R_i^3(z), ~~\widetilde {c}_{i,j}\in \mathbb{C},
\ee
where $\frac{\overline  R_i^3}{\partial z^k}(\Y_1)=0$  for $k=1,\ldots,n$. Moreover, the square  matrix  ${\partial \overline P_i^0\over \partial z^j}(\Y_1)$, $i,j=1,\ldots,r$  is nondegenerate.
 \ep
 \begin{proof}
  The restriction $\overline P_i^0$ of $P_i$ to $Q$ is obtained by setting $z^0=\eta_r$ and $z^k=0$ for $k>n$ in the form \eqref{respoly1}. From the quasihomogeneity of $\overline P_i^0$ and  lemma \ref{first norm}
\be
\overline P_i^0(z^1,\ldots,z^n)=\sum_{a=0}^{\nu_i-1}\sum_{\deg P_i-\deg z^j=a(\eta_r+1)} \widetilde{c}_{i,j,a}(z^r)^{a} z^j+\overline R_i^3(z), ~~\widetilde {c}_{i,j,a}\in \mathbb{C}
\ee
where ${ \overline R_i^3}$ is the  restriction of $R_i^3$ to $Q$. The expressions given in   \eqref{respoly21} imply that ${\partial \overline R^3_i\over \partial z^k}(\Y_1)=0$, $ k=1,\ldots,n$. Note that $\deg P_i-\deg z^j=\nu_i-\eta_j=a(\eta_r+1)$. Using the relation  between the multisets $\exx (\g)$ and $\exx(L_1)$ observed in lemma \ref{relat-expon}, $a$ can only equal $\mu_i$ and  the values of $\eta_j$ are  uniquely determined and depends on $i$. On other words the constants $c_{i,j,a}$ in   \eqref{respoly1}  are nonzero only if $a=\mu_i$. This gives the form \eqref{norm-coord eq}. For the nondegeneracy condition, note that the only possible value for the index $a$  in \eqref{respoly1} is  $a=\mu_i$ and so $z^0$ appear only with the power $\nu_i-\mu_i$. This implies that   ${\partial P_i\over \partial z^j}(\Y_1)={\partial \overline P_i^0\over \partial z^j}(\Y_1)$. Thus  the required matrix is nondegenerate. \end{proof}

Now we give a proof for theorem \ref{nice coordinates} stated in the introduction.

 \begin{proof}[Proof of theorem \ref{nice coordinates}:]
 Writing $\overline P_i^0$ in the form \eqref{arg:sht} and using  the last proposition,  we get  $ \overline P_i^0(q+\lambda K_1)=\overline P_i^0(z^1+ \lambda \delta_{1r},\ldots, z^n+ \lambda \delta_{nr})$ and  $\varrho_i=\mu_i$. We observe that each $\partial_{z^r}^{\mu_i}\overline P_i^0$ is a constant  multiple of $\overline P_i^{\mu_i}$. Hence, the functions  $\partial_{z^r}^{\mu_i}\overline P_i^0$ are Casimirs  of $B_1^Q$ and are in involution with respect to $B_2^Q$ .   Furthermore, $\partial_{z^r}^{\mu_i}\overline P_i^0$  has the form
 \be \partial_{z^r}^{\mu_i}\overline P_i^0=\sum_{\eta_i-\eta_j=\mu_i(\eta_r+1)} \overline{c}_{i,j} z^j+\partial_{z^r}^{\mu_i}\overline R_i^3(z),~\overline{c}_{i,j}\in \mathbb{C},
 \ee
 where $\partial_{z_j}\partial_{z^r}^{\mu_i}\overline R_i^3$ equals 0 at the origin ($z^k=0, \forall k)$. Thus
 \be
 \partial_{z^j}\partial_{z^r}^{\mu_i}\overline P_i^0(0)= {1\over \mu_i!}  {\partial \overline P_i^0\over \partial z^j}(\Y_1), ~ i,j=1,\ldots,r.
 \ee
 We conclude, using  proposition \ref{norm-coord}, that  the matrix   $\partial_{z^j}\partial_{z^r}^{\mu_i}\overline P_i^0$ is nondegenerate. Hence, $ \partial_{z^r}^{\mu_i}\overline P_i^0$ can replace the coordinates $z^i$ on $Q$ for $i=1,\ldots,r$ up to some permutation related to the repetition on $\exx(L_1)$.   Moreover,  using simple linear elimination, we can get the required normalization $t^j= z^j+(\mathrm{non~ linear~ terms})$ where  $t^j$ is  a Casimir of $B_1^Q$. From theorem \ref{rank of B1}, it follows that $t^1,\ldots,t^r$ form a complete set of Casimirs for $B_1^Q$. The fact that $t^1=z^1$
 follows from identifying $t^1 $ with  the Casimir function $\bil Q Q$ and using  $\bil {\ga_1} {L_1}=1$.
 \end{proof}

 We fix the notations $(t^1,\ldots,t^n)$ for the coordinates obtained in theorem \ref{nice coordinates}. Recall that $Z_r(a_s)$  denotes the type of  $L_1$.
 \bc \label{fittness}
 The functions $\overline P_1^0,\ldots,\overline P_{r-s}^0$ are quasihomogeneous polynomials on $t^1,\ldots,t^r$ only.
 \ec
 \begin{proof}
 This follows from the fact that $\mu_i=0$ for $i=1,\ldots,r-s$ and the construction of the coordinates $(t^1,\ldots,t^r)$.
 \end{proof}

\subsection{Integrability and alternative  definitions}

We combine  the theorems stated in section \ref{background1} to get the following useful result.

\bt\label{rank of B1}
The family $\mathbf F$ is complete for every distinguished nilpotent element of semisimple type. In particular, $\rk\,B_1^Q=n-r$ and
 \be\label{def of N}
 N=\{q\in Q: \ker B_1^Q(q)=\ker B_2^Q(q)\}.
 \ee
\et
\begin{proof}
  For regular, subregular and nilpotent elements stated in theorem \ref{mainthm2}, the family $\mathbf F$ is complete  \cite{mypaper5}.   Suppose $L_1$ belongs to the nilpotent orbit  $E_7(a_5)$, $E_8(a_5)$, $E_8(a_6)$, $E_8(a_7)$ or $F_4(a_3)$. We will check that $\rk\, B_\lambda^Q=n-r$ for every $\lambda\in \overline\mathbb C$ and use theorem \ref{bolsinov1}. It is not hard to show that  $\rk\,B_\lambda^Q=n-r$ for $\lambda\in\mathbb C$  \cite{mypaper5}. We need to show that  $\rk\,B_1^Q=n-r$. We verify the equality  by direct computations using proposition \ref{finding rank}  given below. More precisely,  we fixed arbitrary basis $L_i$ for $\g^f$ and $K_1$ such that $L_1+K_1$ is regular semisimple. Then, we found that  the rank of the matrix $\bil {L_i}{[K_1,L_j]}$ equals $n-r$.  The last statement follows from  theorem \ref{bolsinov2}.
\end{proof}

Let us use the special coordinates on $Q$ and denote the entries of the  matrix of the reduced  Poisson structures  by
 \be \label{ldterm}
 F^{ij}_\alpha(t):= \{t^i,t^j\}_\alpha^Q, ~~\alpha=1,2.
 \ee
Then we prove theorem \ref{des of N} stated in the introduction.

 \begin{proof}[Proof of theorem \ref{des of N}:]
 The first definition \eqref{defEqN} of $N$ follows  directly from the structure of the matrices of the Poisson brackets  under the coordinates $(t^1,\ldots,t^n)$. For the second definition \eqref{defEqN1}, we   observe that $dt^1,\ldots,dt^r$ are a basis of  $\ker B_1^Q$ while $d\overline P_1^0,\ldots,d\overline P_r^0$ are basis for $\ker B_2^Q$. However, by construction  $\overline P_1^0,\ldots,\overline P_{r-s}^0$ are polynomials in $t^1,\ldots,t^r$ only. Hence the two kernels coincide exactly on the defined set.

 Now we consider the restriction of the adjoint quotient map  \be
 \Psi^Q(t^1,t^2,\ldots,t^n)=(\overline P_1^0,\ldots,\overline P^0_{r}).
 \ee
 and let $\mathbf J\Psi^Q:={\partial \overline P_i^0 \over \partial t^j} $ denotes its  Jacobian matrix. Then  $N$ is defined by the set of points $t$ where the lower-right $s\times (n-r)$ minor of $\mathbf J\Psi^Q$ is identically 0. From   corollary \ref{fittness}, the upper-right $(r-s) \times (n-r)$ minor  also vanishes by corollary \ref{fittness}. Since,  regular points of $\Psi^Q$ are Zariski dense in $Q$, there exists open dense set $N_0\subseteq N$ such that the left $r\times r$ minor of $\mathbf J\Psi^Q$  is nondegenerate. In particular,  $\overline P_1^0,\ldots,\overline P_r^0$ are independent functions on $N_0$. Hence, $\overline P_1^0,\ldots,\overline P_r^0$ are a part of local coordinates and $\dim N_0\geq r$.  However, the second definition \eqref{defEqN1} of $N$ with corollary \ref{fittness}  implies that $\dim N\leq r$. Thus $\dim N_0=r$ and $(t^1,\ldots,t^r)$ acts as local  coordinates around each point of $N_0$.

 Recall that $B_\lambda^Q$, $\lambda\in \mathbb C$, is of rank $n-r$. Since $N_0$ consists of regular points  the lower-right  $(n-r)\times (n-r)$ minor $F^{\alpha \beta}_\lambda$ of $F^{ij}_\lambda$ is nondegenerate. Thus, Dirac reduction is well defined on $N_0$. However, applying corollary \ref{cor dirac}, the reduced Poisson structure is zero as $t^1,\ldots,t^r$ are in involution with respect to the pencil $B_\lambda^Q$.
 \end{proof}

\section{Algebraic classical $W$-algebra}

In this section, we summarize the construction of  Drinfeld-Sokolov bihamiltonian structure associated to the nilpotent element $L_1$ and $K_1$. Then we will apply Dirac reduction to get a local bihamiltonian structure admitting a dispersionless limit  on the loop space $\N:=\lop N$. This leads to an algebraic classical $W$-algebra on $\N$.

\subsection{Drinfeld-Sokolov reduction}

We consider the loop algebra $\lop \g$ and we extend the bilinear form $\bil . .$ on $\g$ to $\lop \g$ by setting
\begin{equation} (g_1|g_2)=\int_{S^1}\bil {g_1(x)}{g_2(x)} dx;~~~ g_1,g_2 \in \lop \g.
\end{equation}
We use $(.|.)$ to identify $\lop\g$ with $\lop \g^*$.  We define the gradient $\delta \f (g)$ for a functional $\f$ on $\lop\g$
 to be the unique element in
$\lop\g$ satisfying
\begin{equation}
\frac{d}{d\theta}\f(g+\theta
{{\mathrm{w}}})\mid_{\theta=0}=\big(\delta \f(g)|{\mathrm{w}}\big)
~~~\textrm{for all } {\mathrm{w}}\in \lop\g.
\end{equation}
Then,  we introduce standard compatible local Poisson  brackets  $\{.,.\}_1$ and $\{.,.\}_2$ on $\lop \g$ defined for any functionals $\I$ and $\f$ on $\lop \g$ by
 \begin{eqnarray}
\{\f,\I\}_1(g(x)) &:=&\int_{S_1}\bil {[\delta \I(g(x)) ,K_1]}{\delta
\f(g(x))} dx, \\\nonumber
\{\f,\I\}_2(g(x)) &:=& \int_{S^1}\bil {\partial_x\delta \I(g(x))+[\delta \I(g(x)),g(x)]}{ \delta
\f(g(x))} dx.
 \end{eqnarray}
We denote their Poisson structures by  $\mathbb B_1$ and  $\mathbb B_2$, respectively. We mention that  $\mathbb B_2$
  can be interpreted as  the restriction to $\lop \g$ of  Lie-Poisson structure  on the untwisted affine Kac-Moody algebra
 associated to $\g$. In particular, if we expand these Poisson brackets as in   \eqref{bihamileading}, the leading term $\{.,.\}^{[-1]}_1$ is the frozen Lie-Poisson structure $B_1^\g$  and $\{.,.\}^{[-1]}_2$ defines the  Lie-Poisson structure  $B_2^\g$ on $\g$. Moreover, it is easy to show that these Poisson structures  form an exact Poisson pencil  with Liouville   vector field  $\partial_{z^r}$ in the coordinates  defined by \eqref{coord on g}, i.e.,
\be
\{.,.\}_1=\Lie_{\partial_{z^r}} \{.,.\}_2,~~\Lie_{\partial_{z^r}}\{.,.\}_1=0.
\ee

Let us define the affine loop space \be
\Q:=  L_1+\lop {\g^f}.
\ee
Then  Slodowy slice $Q$ is identified with the subspace of constant loops of $\Q$.

\bt \cite{mypaper4}\label{psred}
The space $\Q$ inherits  compatible local Poisson structures   $\mathbb B_2^\Q$ and  $\mathbb B_1^\Q$ from  $\mathbb B_2$ and  $\mathbb B_1$,  respectively. They  can be obtained equivalently by using the bihamiltonian reduction with Poisson tensor procedure, Dirac reduction and the generalized  Drinfeld-Sokolov reduction. Moreover, the leading terms of the bihamiltonian structure on $\Q$ can be identified with the bihamiltonian structure $B_2^Q$ and $B_1^Q$ on $Q$.
\et

Details on bihamiltonian reduction can be found in \cite{CMP}.  Drinfeld-Sokolov reduction is initiated and applied for regular nilpotent elements in \cite{DS}. Generalizations to other nilpotent elements is obtained in   \cite{gDSh2},\cite{BalFeh1} (see also \cite{mypaper}). The relation between Drinfeld-Sokolov reduction and bihamiltonian reduction in the case of regular nilpotent elements is treated in  \cite{CP} and \cite{Pedroni2}. In \cite{CP}, the Poisson tensor procedure is also initiated (also called the method of transverse subspace in \cite{LPM}). The relation between Drinfeld-Sokolov reduction and Dirac reduction is also proved in \cite{BalFeh1}. See \cite{kac} and references therein, for more recent  development and tools used to study  Drinfeld-Sokolov reduction.

We let  $\{.,.\}_1^\Q$ and $\{.,.\}_2^\Q$ denote the Poisson brackets  defined by $\mathbb B_1^\Q$ and  $\mathbb B_2^\Q$, respectively.

In what follows, we review  Drinfeld-Sokolov reduction. We identify $\lop \g$ with the space of operators of the form $\partial_x+ g$, $g\in \lop \g$, and $\Q$ with the subspace of operators of the form $\partial_x+q+ L_1$, $q\in \lop {\g^f}$. Let $\B$ denote the subspace of  operators of the form
\begin{equation}\label{op:DS}
\L=\px +b+L_1\qquad\textrm{where
} b\in \lop\bneg,~~\bneg:=\bigoplus_{i\leq 0} \g_i.
\end{equation}
There is a natural action of the adjoint group of
$\lop\nneg$, $\nneg:=\bigoplus_{i<0}\g_i$,  on $\B$  defined by
\begin{equation}\label{gauge fix}
 ({\mathrm{w}},\L)\to (\exp{\ad {\mathrm{w}}})\, \L \textrm{ for all }
{\mathrm{w}}\in \lop\nneg \textrm{ and } \L\in \B.
\end{equation}
Moreover, for any operator $\L\in \B$ there is a unique element ${\mathrm{w}} \in
\lop\nneg$ such that
\begin{equation}\label{op:fixing} \lcan:=\px +q+ L_1=(\exp{\ad {\mathrm{w}}}) \L
\end{equation}
where $ q\in \lop{\g^f}$. Hence, $q$ and ${\mathrm{w}}$ are differential polynomials in the coordinates of $b$.  The entries of $q$ give a set of  generators
of the ring $R$ of differential polynomials invariant under the action \eqref{gauge fix} . More precisely, if we write
 \begin{equation}
b=\sum_{i=1}^n\sum_{I=0}^{\eta_i} b_I^i (x){1\over I!} \ad_{L_1}^{I} \ga_i,~~q=\sum_{i=1}^n z^{i}(x)  \ga_i ~~\textrm{and } {\mathrm{w}}=\sum_{i=1}^n\sum_{I=1}^{\eta_i} {\mathrm{w}}_I^i (x){1\over I!} \ad_{L_1}^{I} \ga_i,
\end{equation}
then  equation \eqref{op:fixing} reads
\be \label{op:fixing1}
q-[{\mathrm{w}},L_1]=b-{\mathrm{w}}_x+[{\mathrm{w}},b]+\sum_{i>0} {1\over i+1!}\ad_{\mathrm{w}}^i (-{\mathrm{w}}_x+[{\mathrm{w}},b]+[{\mathrm{w}},L_1]).
\ee
Using Dynkin grading and the fact that $\g^f\oplus [\nneg,L_1]=\bneg$, we get recursive equations defining
 the coordinates of $q$ as differential polynomials on the  coordinates of $b$.  Moreover, if we assign degree
$\partial_x^k b_{J}^{i}$ equals $k+\eta_i-J+1 $, then $z^i(x)$ is a quasihomogeneous polynomial of degree $\eta_i+1$.
The  set of
functionals $\mathcal{R} $ on $\Q$ are the functionals on $\B$ with
densities belonging to the ring $R$. It follows that $\mathcal{R}$ is closed Poisson subalgebra with respect to the Poisson brackets $\{.,.\}_2$ and $\{.,.\}_1$. Thus, the reduced Poisson pencil $\{.,.\}^\Q_\lambda:=\{.,.\}_2^\Q+\lambda\{.,.\}_1^\Q $ can be obtained by apply the Leibniz rule
\begin{equation}\label{leibniz rule}
\{z^u(x),z^v(y)\}_\lambda^\Q:={\partial z^u(x)\over \partial (b^{i}_I)^{(k)}}\partial_x^k\Big({\partial z^v(y)\over \partial (b_{J}^{j})^{(l)}} \partial_y^n\big(\{b_{I}^{i}(x),b_{J}^j(y)\}_\lambda\big)\Big)
\end{equation}
where
\begin{eqnarray}
\{b_{I}^{i}(x),b_{J}^{j}(y)\}_\lambda &=&{1\over \Theta_{I}^{i}}{1\over \Theta_{J}^{j}} \Big( \bil  {\ad_f^J L_j}{\ad_f^I L_i}\partial_x+\bil {b(x)+\lambda K_1} {[\ad_f^J L_j,\ad_f^I L_i]}\Big)\delta(x-y)\\ \nonumber
&=&{1\over \Theta_{I}^{i}}{1\over \Theta_{J}^{j}} \Big( \bil  {\ad_f^J L_j}{\ad_f^I L_i}\partial_x+ \bil {\ad_f^I L_i} {[b+\lambda K_1,\ad_f^J L_j]}\Big)\delta(x-y)
\end{eqnarray}
and $\Theta_{J}^{j}:=(-1)^J {{\eta_j}\choose{J}}$. We  will use these formulas in next sections to analyse the leading terms of  $\mathbb B_2^Q$ and $\mathbb B_1^Q$.

We end this section by finding the linear terms of the generators of the invariant  ring $R$.

\bp \label{linear part}  The  linear terms of each $z^i(x)$ equals
\be \sum_{I=0}^{\eta_i} {(-1)^I\over I!} \partial_x^I b_{I}^{i}.\ee
In particular, $z^r(x)$ is the only generator of $R$ depends on $b_0^r(x)$  and this dependence is linear. Moreover,  all $z^i(x)$ do not depend on derivatives of $b_0^r(x)$.
\ep
\begin{proof}
The second part of the statement follows from the quasihomogeneity of the generators $z^i(x)$ of $R$. To find linear terms of each $z^i$, we introduce spectral parameter $\epsilon$ and set $\L(\epsilon)=\partial_x+\epsilon b+ L_1$. Let ${\mathrm{w}}(\epsilon)$ and $\lcan(\epsilon)$ be the corresponding operators. Then $\L(0)=\partial_x+L_1$, ${\mathrm{w}}(0)=0$ and $\lcan(0)=\L(0)$. Therefore, differentiating the relation
\be
\lcan(\epsilon)=\L(\epsilon)+ [\fn(\epsilon),\L(\epsilon)]+{1\over 2}[\fn(\epsilon),[\fn(\epsilon),\L(\epsilon)]]+\ldots
\ee
with respect to $\epsilon$ and evaluating at $\epsilon=0$ we get
\begin{eqnarray}
q'(0)&=&b+[{\mathrm{w}}(0), \L'(0)]+ [{\mathrm{w}}'(0),\partial_x+L_1]\\\nn
&=& b+[{\mathrm{w}}'(0),\partial_x+L_1]\\\nn
&=& b -{\mathrm{w}}_x'(0)+ [{\mathrm{w}}'(0),L_1].
\end{eqnarray}
Note that $[{\mathrm{w}}'(0),L_1]$ does not contribute to $q'(0)$.  Then  the coordinate  of $\ga_i$ gives
\be
(z^i)'(0)=b_{0}^{i}-({\mathrm{w}}_x'(0))_{0}^{i}
\ee
where we write ${\mathrm{w}}'(0)=\sum_{i=1}^n\sum_{I>0} ({\mathrm{w}}'(0))_I^i{1\over I!} \ad_{L_1}^{I} \ga_i$. Then the coefficients of ${1\over I!}ad_{L_1}^I\ga_i $ for $I>0$ give the recursive relations
\be [({\mathrm{w}}'(0))_{I-1}^{i},L_1]-({\mathrm{w}}_x'(0))_{I}^i+b_{I}^i=0
\ee
which leads to
\be
 ({\mathrm{w}}'(0))_{I-1}^{i}={1\over I+1}(-({\mathrm{w}}_x'(0))_{I}^{i}+b_{I}^{i}).
 \ee
For example
\begin{eqnarray}
 ({\mathrm{w}}'(0))_{\eta_i-1}^{i}&=&{1\over \eta_i}(b_{\eta_i}^{i}),\\ \nn
  ({\mathrm{w}}'(0))_{\eta_i-2}^{i}&=&{1\over \eta_i-1}(-{1\over \eta_i}(\partial_x b_{\eta_i}^{i})+b_{\eta_i-1}^{i}).
\end{eqnarray}
These recursive relations lead to
\be
(z^i)'(0)=\sum_{I=0}^{\eta_i} {(-1)^I\over I!} \partial_x^I b_{I}^{i}.
\ee
\end{proof}

Recall that the coordinates $(t^1,\ldots,t^n)$ of $Q$ developed in theorem \ref{nice coordinates} are  quasihomogeneous polynomials in the coordinates $(z^1,\ldots,z^n)$. Thus we get the following corollary by construction.

\bc \label{linear for t} Proposition \ref{linear part} is valid when we replace $z^i(x)$ by $t^i(x)$.
\ec

\subsection{Further reduction}

In this section, we reduce Drinfeld-Sokolov bihamiltonian structure to $\N$  and analyze the leading term  using the coordinates  $(t^1,\ldots,t^n)$ obtained by theorem  \ref{nice coordinates}.

\bp\label{exactness}
The reduced bihamiltonian structure on $\Q$ is exact with Liouville vector field $ \partial_{t^r}$. The Poisson bracket with $t^1$ preserve the relations  defining classical $W$-algebra, i.e.,
\begin{eqnarray}\label{leading terms1}
\{t^1(x),t^1(y)\}_2^Q&=&  c \delta^{'''}(x-y) +2 t^1(x) \delta'(x-y)+ t^1_x\delta(x-y), \\\nonumber
\{t^1(x),t^i(y)\}_2^Q &=& (\eta_i+1) t^i(x) \delta'(x-y)+ \eta_i t^i_x \delta(x-y), ~i=2,\ldots,n.
\end{eqnarray}
for some nonzero constant $c$.
\ep

\begin{proof}
We take  $t^1(z),\ldots, t^n(z)$ as generators for the invariant ring $R$. By corollary \ref{linear for t},   $t^r(x)$  is the only invariant which depends on $b_{0}^{r}(x)$. This  implies that the invariant  $t^{r}(x)$  appears in the expression of $  \{t^i(x),t^j(y)\}^{\Q}_2$ only if, when using the Leibniz rule \eqref{leibniz rule}, we encounter  terms  of $\{.,.\}_2$ depend explicitly on  $b_{0}^{r}(x)$. Thus $  \{t^i(x),t^j(y)\}^{\Q}_2$ is at most linear on $z^r(x)$ and its derivatives. But the bihamiltonian structure on $\lop \g$ is exact and $\{.,.\}_1$ is obtained from $\{.,.\}_2$  by the shift along $b_0^r$. Hence, $  \{t^i(x),t^j(y)\}^{\Q}_1$ is obtained by the shift of $  \{t^i(x),t^j(y)\}^{\Q}_2$ along $t^r(x)$, i.e., substituting $t^r(x)$ by $t^r(x)+\epsilon$ and evaluate  ${d\over d\epsilon}|_{\epsilon=0}$. Therefore, $  \{.,.\}^{\Q}_1$ does not depend on $t^r(x)$ or its derivatives. From the work in  \cite{BalFeh1}, the  reduced Poisson bracket $\{.,.\}_2^\Q$ is a classical $W$-algebra in the coordinates $(z^1,\ldots,z^n)$, i.e., it satisfies the identities \ref{walgebra2}.
Then the argument for identities \eqref{leading terms1}  will be similar to the one given in the proof  of proposition \ref{flat in s} below.
\end{proof}

Then  theorem  \ref{classical in N} gives  compatible local Poisson brackets $\{.,.\}_\alpha^\N$, $\alpha=1,2$ on the loop space $\N=\lop N$ of the space of common equilibrium points $N$. The proof is as follows.

\begin{proof}[Proof of theorem \ref{classical in N}]
From theorem  \ref{des of N},   the leading terms of  $\{.,.\}_k^\Q$, $k=1,2$ have the form  \be \{t^i(x),t^j(y)\}^{[-1]}_k = F^{ij}_k(t(x))\delta(x-y),~\ee
where  $F^{i\alpha}_1(t)=0$  and $N$ is defined by  $F^{i\alpha}_2(t)=0$, $1\leq i\leq r$ and $r+1\leq \alpha\leq n$. Thus, $\{.,.\}_\lambda^\Q$ satisfies the hypothesis of proposition \ref{dirac fromula} with the coordinates $(t^1,\ldots,t^r)$ on $N$. Using corollary \ref{cor dirac}, the reduced local Poisson bracket $\{.,.\}_\lambda^\N$ on $\N$ is obtained by setting $\{t^i(x),t^j(y)\}_\lambda^\N$ equals  $\{t^i(x),t^j(y)\}_\lambda^\Q$ and substitute   the variables $t^{i}$, $i>r$ by solutions of the polynomial equations  $F_2^{i\alpha}=0$ defining $N$. In particular, $\{t^i(x),t^j(y)\}_\lambda^\N$ is an algebraic local Poisson bracket and it is linear in $\lambda$. This leads to compatible local Poisson brackets $\{.,.\}_2^\N$ and $\{.,.\}_1^\N$ on $N$ where the former still satisfies  the identities \eqref{leading terms1} defining classical $W$-algebras. From theorem  \ref{des of N} again, they both admit a dispersionless limit. Note that the defining equation  $F_2^{i\alpha}=0$ of $N$ do not depends on $t^r$. Thus, from proposition \ref{exactness}, the reduced Poisson brackets  form an exact Poisson pencil.    \end{proof}

As in the introduction, we write the leading terms of $\{.,.\}^{\N}_{\alpha}$, $\alpha=1,2$, in the form
\be\label{reducedPB notations}
  \{t^u(x),t^v(y)\}^{[0]}_\alpha = \O^{uv}_\alpha(t(x)) \delta' (x-y)+ \Gamma_{\alpha k}^{uv}(t(x)) t_x^k \delta (x-y), ~1\leq u,v\leq r.\ee
In the remainder of this section, we want to prove that   the  determinate of the matrix $\O^{uv}_1(t)$  is nonzero constant. For this end, we write  \be [K_1, \ad_f^J L_j]=\sum_t \Delta_{j}^{Jt} {1\over T!}\ad_{L_1}^T \ga_t;~~ T=\eta_t+\eta_j-J-{\eta_r}\geq 0\ee
where $T $ is constrained  by the Dynkin grading of $\g$.
Then  the values of $\{.,.\}_1$ on the coordinates of $b$ are given by
\be  \{b_{I}^{i}(x),b_{J}^{j}(y)\}_1= {1\over \Theta_{J}^{j} } \delta^{IT}\delta^{it} \Delta_{j}^{Jt}\delta(x-y).\ee
Thus, we get the following formula for the brackets
\be \label{ijcontstrain}
\{b_{I}^{i}(x),b_{J}^{j}(y)\}_1=  {\Delta_{j}^{Ji} \over \Theta_{J}^{j} }  \delta(x-y),~~~I=\eta_i+\eta_j-J-{\eta_r}\ee
  where $\Delta_{j}^{Jt}$ possibly equals 0.
  Expanding using the Leibniz rule, we get
  {\small \begin{eqnarray}\label{general exp}
\{t^u(x),t^v(y)\}_1^\Q&=& \sum_{i,j}\sum_{l,h}{\Delta_{j}^{Ji} \over \Theta_{J}^{j} }{\partial t^u(x) \over \partial (b_{I}^{i})^{(l)} }\partial_x^l\Big({\partial t^v(y) \over \partial (b_{J}^{j})^{(h)} }\partial_y^h\delta(x-y)\Big),~I=\eta_i+\eta_j-J-{\eta_r}\\\nonumber
&=&  \sum_{i,j}\sum_{l,h,\alpha,\beta} (-1)^h{h\choose \alpha} {l\choose \beta}{\Delta_{j}^{Ji} \over \Theta_{J}^{j} }{{\partial t^u(x) \over \partial (b_{I}^{i})^{(l)} }\Big({\partial t^v(x) \over \partial (b_{J}^{j})^{(h)} }\Big)^{(\alpha+\beta)}\delta^{(h+l-\alpha-\beta)}(x-y)}.
\end{eqnarray}}
Here we omitted the ranges of the indices since no confusion can arise. We observe that the value of $\O^{uv}$ is contained in the expression
\begin{eqnarray}\label{g1:matrix}
\f^{uv}_1&=&\sum_{i,J}\sum_{h,l}(-1)^h (l+h) {\Delta_{j}^{Ji} \over \Theta_{J}^{j} }{\partial t^u(x) \over \partial (b_{I}^{i})^{(l)} }\Big({\partial t^v(x) \over \partial (b_{J}^{j})^{(h)} }\Big)^{h+l-1},~~I=\eta_i+\eta_j-J-{\eta_r}
\end{eqnarray}

 \bl   \label{homg of g1}  The matrix $\O^{uv}_1(t)$ is lower antidiagonal with respect to $\exx(L_1)$ and the antidiagonal entries are constants. In other words, $\O^{uv}_1(t)$  is constant  if $\eta_u+\eta_v= \eta_r+1$
 and equals zero if $\eta_u+\eta_v<\eta_r+1$.
\el
\begin{proof}
Assume   $t^u(x)$ and $t^v(x)$ are  quasihomogeneous of  degree $ \eta_u+1$ and $\eta_v+1$, respectively. Then $\f^{uv}_1$ is a quasihomogeneous polynomial of degree \[ \eta_{u}+1+\eta_v+1-(\eta_i-I+l+1)-(\eta_j-J+h+1)+h+l-1=\eta_u+\eta_v-{\eta_r}-1\]
\end{proof}

Recall  that from  the construction of the  coordinates $(t^1,\ldots,t^r)$ and the second part of proposition \ref{opp Cartan prop},  the entry $\O^{uv}_1$ in case $u+v=r+1$ implies that $\eta_u+\eta_v=\eta_r+1$ and $\bil {\Y_u}{\Y_v}=\eta_r+1$

\bp\label{nondeg}
The antidiagonal entries  of  $\O^{uv}_1$ with respect to the set $\exx(L_1)$ equal  $\eta_r+1$ in case $u+v=r+1$ and zero otherwise. In particular,
 $\O^{uv}_1$ is nondegenerate and its determinant equals to $(\eta_r+1)^r$.
\ep
\begin{proof}
We need only to examine the entry  $\O^{uv}$ where $t^u$ and $t^v$ are quasihomogeneous of degree $\eta_u+1$ and ${\eta_r}-\eta_u+2$, respectively.  The expression \eqref{g1:matrix} yields the constrains
\begin{eqnarray}
\eta_i+1-I\leq \eta_u+1 &\Rightarrow& J\leq \eta_u+\eta_j-{\eta_r} \\\nn
\eta_j+1-J\leq {\eta_r}-\eta_u+2&\Rightarrow & \eta_j+\eta_u-{\eta_r}-1\leq J.
\end{eqnarray}
 Hence  $J$ equals  $\eta_u+\eta_j-{\eta_r} -1$ or $\eta_u+\eta_j-{\eta_r}$. Consider  $J=\eta_u+\eta_j-{\eta_r} -1$. Then  $\deg(b_J^j)^{(h)}=\eta_j-J+1+h=\deg t^v+h$. This forces $h=0$ and $t^v$ is linear in $b_J^j$. Therefore, from proposition \ref{linear part}, $j=v$ and $J=0$ which leads to ${\partial t^v(x) \over \partial (b_{J}^{j})^{(h)} }=1$. Also
\be \deg (b_{I}^{i})^{(l)}=\eta_i-I+h+1=\eta_i-(\eta_i+\eta_j-J-\eta_r)+l+1=\eta_u+l=\deg t^u.\ee
Thus the only possible value for $l$ is 1. Note that $I=\eta_i-\eta_u+1$. Hence, $\deg t^u= \deg (b_I^i)'$ and $t^u$ is linear in $(b_I^i)'$. Then $i=u$ and $I=1$ and from proposition \ref{linear part},  ${\partial t^u(x) \over \partial (b_{I}^{i})^{(l)} }=-1$. Therefore,  the case   $ J=\eta_u+\eta_j-{\eta_r} -1$, the expression  \eqref{g1:matrix} contribute to $\O^{uv}_1$ with the value
$-{\Delta_v^{0u}\over \Theta^{v}_0}=-\Delta_v^{0u}$ since  $J=\eta_u+\eta_v-{\eta_r}-1=0$. By  definition,
\be
-\Delta_v^{0u}={1\over \eta_u}\bil {ad_f L_u}{[K_1,L_v]}=\eta_v \delta_{u+v,r+1}.
\ee
A similar analysis when  $J=\eta_u+\eta_j-{\eta_r}$  leads to the value $\eta_u \delta_{u+v,r+1}$. By the normalization of $\Y_i$, it follows that the value of $\O^{uv}_1$   equals  $\eta_u+\eta_v={\eta_r}+1$ when  $u+v=r+1$ and zero otherwise. The determinant   of the matrix $\O^{uv}_1$  follows accordingly.
\end{proof}
\bc \label{nondeg222}
The  matrix $\O_2^{uv}(t)$ is nondegenerate  on $\N$.
\ec
\begin{proof}
It follows from the exactness of the Poisson pencil, i.e.,   $ \O^{uv}_1(t)=\partial_{t^{r}} \O_2^{uv}(t)
$.
\end{proof}

Recall  the duality  of the multiset $\overline \exx(L_1)$ stated on corollary \ref{Ecompl}. Then the following proposition  is useful  to find the rank of $ B_1^Q$. Note that the proof depends only on the linear part of the invariants $t^i(x)$.
\bp \label{finding rank}
The  matrix $F^{uv}_1(t)$, $u,v=1,\ldots,n$ is a lower antidiagonal in the sense that  $F^{uv}_1(t)=0$ if $\eta_u+\eta_v<\eta_r$. In particular, if $\eta_u+\eta_v=\eta_r$ then \be F^{uv}_1(t)=\bil {L_{u}}{[K_1,L_v]},\ee
and if $\eta_u+\eta_v=\eta_r+1$ then $F^{uv}_1(t)=0$
\ep
\begin{proof}
Note that the value of the matrix $F^{uv}_1(t)$ is contained in the expression
\be\label{abexpr}
\sum_{i,J}\sum_{h,l} (-1)^h {\Delta_{j}^{Ji} \over \Theta_{J}^{j} }{\partial t^u(x) \over \partial (b_{I}^{i})^{(l)} }\Big({\partial t^v(x) \over \partial (b_{J}^{j})^{(h)} }\Big)^{h+l},~~I=\eta_i+\eta_j-J-{\eta_r}.
\ee
Then the proof will be   similar  to the proof of  lemma \ref{homg of g1} and  proposition \ref{nondeg}. The degree of this expression is $\eta_u+\eta_v-\eta_r$. Thus the matrix will be lower antidiagonal as claimed. Let us assume $\eta_v+\eta_u=\eta_r$. Then  the only possible value for $J$ is $\eta_u+\eta_j-\eta_r$. We also find   $h$ (resp. $j$, $l$ and $i$) must equal $0$ (resp. $v$, $0$ and $u$). Therefore, $J=0$ and the  expression \eqref{abexpr} will be $\Delta^{0 u}_v=\bil {L_{u}}{[K_1,L_{v}]}$. For the last statement, note that $F^{uv}_1(t)$ is a polynomial \cite{mypaper4} and there is no variable of degree 1.
\end{proof}


\section{Algebraic Frobenius manifold} \label{Alg Frob mani}
In this section, we obtain the promised algebraic Frobenius structure and give examples

\subsection{General construction}

We consider the flat pencil of metrics on $N$ consists of $\O^{uv}_1(t)$ and $\O^{uv}_2(t)$ which is afforded by theorem \ref{DN thm}, theorem \ref{classical in N}, proposition \ref{nondeg} and corollary  \ref{nondeg222}. From the exactness of Poisson pencil on $\N$ and defining equations of $W$-algebra given in proposition  \ref{exactness}, we have
\be \Lie_{\partial_{t^r}}\O^{uv}_2=\O^{uv}_1,~~~g^{1u}_2(t)=(\eta_u+1)t,  ~ \Gamma^{1j}_{2 k}(t)=\eta_j \delta^j_k.
\ee
Recall that we assign degree $t^u$ equals $\eta_u+1$.

\bp \label{homg of g2}
Each entry  $\O^{uv}_2(t)$ is quasihomogeneous of degree $\eta_u+\eta_v$ while  ${\Gamma}^{uv}_{2k}(t)$ is quasihomogeneous of degree $\eta_u+\eta_v-(\eta_k+1)$.
\ep
\begin{proof}
First part follows  from the proof of lemma \ref{homg of g1}. Analysing the coefficient of $\delta(x-y)$ is the expression \eqref{general exp} leads to the degree of ${\Gamma}^{uv}_{2k}(t)$.
\end{proof}

\bp \label{flat in s}
There exist    a quasihomogeneous polynomial change of coordinates of the form
\be  \label{ffflattt} s^i=t^i+\mathrm{non~ linear~ terms}\ee  such that  the matrix $\O_1^{uv}(s)=(\eta_r+1)\delta^{u+v,r+1}$. Furthermore, in these coordinates the metric $\O^{uv}_2(s)$ and its Christoffel symbols  preserve the identities
\begin{equation}
\O^{1,v}_2(s)= (\eta_v+1) s^v, ~ \Gamma^{1v}_{2 k}(t)=\eta_v \delta^v_k.
\end{equation}
\ep
\begin{proof}
A local flat coordinates of the metric $\O_1^{uv}(s)$ exist at each point of $N$  and can be found by solving the system \cite{DCG}
\be
\O_1^{uv} \partial_{t^u}\partial_{t^k} s +\Gamma^{uv}_{1 k}\partial_{t^v} s=0, ~u,k=1,\ldots,r.
\ee
First, we  search for a quasihomogeneous change of coordinates in the form  $s^i=s^i(t^1,\ldots,t^r)$ with $\deg s^i=\deg t^i$ such that  the matrix $\O^{uv}_1(s)$ is constant antidiagonal with respect to the set $\exx(L_1)$. The proof of its existence  can be obtained by following the proof of a similar statement in  (\cite{DCG}, corollary 2.4).  Note that we can write $s^i$ in the form \eqref{ffflattt} using eliminations. But then, after reordering, we can apply proposition \ref{nondeg} to get $\O_1^{uv}(s)=(\eta_r+1)\delta^{u+v,r+1}$.  For the second part of the statement,  we need only to show that  \begin{equation}
\O^{1,i}_2(s)= (\eta_i+1) s^i, ~ \Gamma^{1j}_{2k}(s)=\eta_j \delta^j_k.
\end{equation}
Let us  introduce  the Euler vector field
\beq
E':=\sum_i (\eta_i+1) t^i { \partial_{t^i}}.
\eeq
Then the formula for change of coordinates   gives
\beq
\O^{1j}_2(s)={\partial_{t^a} s^1 } {\partial_{t^b} s^j}~ \O_2^{a b}(t)= E'(s^j)=(\eta_j+1) s^j.
\eeq
Here the last equality comes from quasihomogeneity of the coordinates $s^i$. For    $\Gamma^{1j}_{2k}(t)$, the change of coordinates has the following formula
\be
\Gamma^{ij}_{2k}(s) d s^k=\Big({\partial_{t^a} s^i } {\partial_{t^c}\partial_{t^b} s^j} \O_2^{a b}(t)+ {\partial_{t^a} s^i} {\partial_{t^b} s^j  } \Gamma^{a b}_{c}(t)\Big) d s^c.
\ee
But then we get
\begin{eqnarray}
\Gamma^{1j}_{2k} d s^k&=&\Big( E' ({\partial_{t^c} s^j  })+  {\partial_{t^b} s^j  } \Gamma^{1 b}_{2c}\Big) d t^c\\\nonumber
&=& \Big( (\eta_j-\eta_c){\partial_{t^c} s^j  }+  \eta_c {\partial_{t^c} s^j  } \Big) d t^c=\eta_j {\partial_{t^c} s^j } d t^c =\eta_j d s^j.
\end{eqnarray}

\end{proof}

From proposition \ref{flat in s}, we  can assume without loss of generality that the  coordinates $t^i$ are the flat coordinates for $ \O^{ij}_1$. Then we get a regular quasihomoegenius  flat pencil of metrics of degree $\frac{\eta_r-1}{\eta_r+1}$ formed by  $ \O^{ij}_1$ and $ \O^{ij}_2$ on $N$ as theorem \ref{my thm} states.

\begin{proof}[Proof of theorem \ref{my thm}]
In the notation of equations (\ref{tau flat pencil}), we set $\tau :={1\over {\eta_r}+1}t^1$. Then
\begin{eqnarray}\label{Euler of alg}
 E&:=& \O^{ij}_2 {\partial_{t^j} \tau  }~{\partial_{t^i}  }={1\over {\eta_r}+1} \sum_{i} (\eta_i+1) t^i{\partial_{t^i} },\\\nonumber
 e &:=&  \O^{ij}_1 {\partial_{t^j} \tau }~{\partial_{t^i}  }={\partial_{t^{r}} }.
  \end{eqnarray}
The identities  $[e,E]=e$,
$\Lie_{\partial_{t^r}}\O^{uv}_2=\O^{uv}_1$ and $\Lie_{\partial_{t^r}}\O^{uv}_1=0$ are fulfilled. We also obtain from proposition \ref{homg of g2} that
  \be
\Lie_E \O^{ij}= E(\O^{ij}_2)-{\eta_i+1\over {\eta_r}+1}\O^{ij}_2-{\eta_j+1\over {\eta_r}+1} \O^{ij}_2={-2\over {\eta_r}+1} \O^{ij}_2=(d-1) \O^{ij}.
\ee
 We also have the regularity condition since the (1,1)-tensor $R_i^j$  has the entries
\begin{equation}
  R_i^j = {d-1\over 2}\delta_i^j + {\nabla_1}_iE^j = {\eta_i \over {\eta_r}+1} \delta_i^j.
\end{equation}
\end{proof}

Now we can prove the main result,  theorem \ref{main thm}.

\begin{proof}[Proof of theorem \ref{main thm}]
It follows from  theorems \ref{my thm} and \ref{dub flat pencil}  that $N$  has a natural Frobenius structure of charge  ${\eta_r}-1\over {\eta_r}+1$. This Frobenius structure is algebraic   since  the potential $\mathbb{F}$ is constructed using   equations \eqref{frob eqs} and from theorem \ref{classical in N}  the matrix $\O^{uv}_2$ may contain variables $t^k,k>r$ which are solution  of the polynomial equations \eqref{defEqN} defining $N$. The Euler vector field is given by the formula \eqref{Euler of alg}. By construction, different choices of a representative $L_1$ or transverse subspace other than Slodowy slice will lead to the same   Frobenius structure.
\end{proof}

\subsection{Examples}

\subsubsection{Regular nilpotent orbits}
Suppose $L_1$ is a regular nilpotent element in $\g$. Then the multisets $\exx(L_1)$ and $\exx(\g)$ coincide. In this case, we get the standard   Drinfeld-Sokolov reduction \cite{DS} on  Slodowy slice $\Q$ and the local bihamiltonian structure  admits a dispersionless limit. Hence, the space of common equilibrium points $N$ equals  $Q$. The algebraic Frobenius manifold is polynomial. It coincides \cite{DLZ} with the polynomial Frobenius manifold constructed by Dubrovin on the orbit spaces of the underlined Weyl group \cite{DCG}. The  construction using the methods of this article was also obtained in \cite{mypaper1}.

\subsubsection{Subregular nilpotent orbits}

A nilpotent elements is called subregular if $\dim \g_0=r+2$. The set of all subregular nilpotent elements form one  nilpotent orbit which exists in any complex simple Lie algebra.  However, not all subregular nilpotent elements of simple Lie algebras are of semisimple type. Which was wrongly assumed in the article $\cite{mypaper3}$. Only the subregular nilpotent elements of type $D_4(a_1)$, $F_4(a_1)$, $E_6(a_1)$, $E_7(a_1)$ and $E_8(a_1)$ are of semisimple type. Hence, all statements in  \cite{mypaper3} are valid only when  considering those cases.
Let $L_1$ be a subregular nilpotent element of semisimple type. Then Slodowy slice $Q$ is of dimension $r+2$. In \cite{mypaper3}, the set of  common equilibrium points  $N$  was  defined in terms of the invariant polynomials $P_1,\ldots,P_r$ using the normalization of the transverse Lie-Poisson bracket $\{.,.\}^Q_2$ obtained  in \cite{DamSab}.  Moreover, the article \cite{mypaper3}  contains in detail the construction of the potential of the algebraic Frobenius manifold associated to $D_4(a_1)$. So we are not keen to repeat writing this example here.  We also constructed the potential associated with $E_8(a_1)$, but it results in a huge polynomial in $8$ variables (consist of 303 monomials) with vast numbers and by all means unpublishable \cite{myarixE8}.  A simpler formula for this potential  appears in \cite{DinarSek}.


\subsubsection{Nilpotent element of type  $F_4(a_2)$}
We use minimal representation of $F_4$ which  is given by  square matrices of size 27. The following computations can be verified  using any computer algebra systems. Below $\epsilon_{i,j}$ denote the standard basis of the set of square matrices of size 27. To simplify the notation we use $E_{c_1c_2c_3c_4}$ to denote the root vector corresponding to the root ${c_1\al_1+c_2\al_2+c_3\al_3+c_4\al_4}$ while $F_{c_1c_2c_3c_4}$ for the root vector corresponding to  the  negative root.  We always set  $F_{c_1c_2c_3c_4}$ equals the transpose of the matrix $E_{c_1c_2c_3c_4}$. Then the simple root vectors are
\begin{eqnarray}
E_{0001} & := & -\epsilon _{4,5}+\epsilon _{7,8}+\epsilon _{9,11}+\epsilon _{20,22}+\epsilon _{21,6}+\epsilon _{23,24},\\\nonumber
E_{0010}& := & -\epsilon _{3,4}+\epsilon _{8,10}+\epsilon _{11,13}+\epsilon _{18,20}+\epsilon _{19,21}+\epsilon _{24,25}, \\\nonumber
E_{0100} & := & -\epsilon _{2,3}-\epsilon _{4,7}+\epsilon _{5,8}+\epsilon _{6,24}+\epsilon _{10,12}+\epsilon _{13,15}+\epsilon _{13,16}+\epsilon _{15,18}+\epsilon _{16,18}+\epsilon
   _{17,19}+\epsilon _{21,23}+\epsilon _{25,26},\\\nonumber
E_{1000} & := &  -\epsilon _{1,2}-\epsilon _{7,9}-\epsilon _{8,11}-\epsilon _{10,13}+\epsilon _{12,14}-\epsilon _{12,15}-\epsilon _{14,17}+\epsilon _{15,17}+\epsilon
   _{18,19}+\epsilon _{20,21}+\epsilon _{22,6}+\epsilon _{26,27}.\\\nonumber
\end{eqnarray}
We construct the remaining root vectors by setting
\begin{center}
\begin{tabular}{|c|c|c|}\hline
 $ E_{0011}=[E_{0001},E_{0010}]$   &  $ E_{0110}=[E_{0010},E_{0100}]$   &   $ E_{1100}=[E_{0100},E_{1000}]$   \\ & & \\\hline
  $ E_{0111}=[E_{0011},E_{0100}]$   &  $ E_{0210}=[E_{0100},E_{0110}]$   &   $ E_{1110}=[E_{1000},E_{0110}]$   \\ & & \\\hline
   $ E_{0211}=[E_{0111},E_{0100}]$   &  $ E_{1111}=[E_{1110},E_{0001}]$   &   $ E_{1210}=[E_{1110},E_{0100}]$   \\ & & \\\hline
  $ E_{0221}=[E_{0211},E_{0010}]$   &  $ E_{1211}=[E_{1111},E_{0100}]$   &   $ E_{2210}=[E_{1210},E_{1000}]$   \\ & & \\\hline
   $ E_{1221}=[E_{0221},E_{1000}]$   &  $ E_{2211}=[E_{1211},E_{1000}]$   &   $ E_{1321}=[E_{1221},E_{0100}]$   \\ & & \\\hline
  $ E_{2221}=[E_{2211},E_{0010}]$   &  $ E_{2321}=[E_{2221},E_{0100}]$   &   $ E_{2421}=[E_{2321},E_{0100}]$   \\ & & \\\hline
  $ E_{2431}=[E_{2421},E_{0010}]$   &  $ E_{2432}=[E_{2431},E_{0001}]$   &     \\ & & \\\hline

\end{tabular}
\end{center}

We fix the following $sl_2$-triple, where the nilpotent element $L_1$ is of type $F_4(a_2)$

\begin{eqnarray}
L_1&=& E_{0010}+E_{0011}+E_{0110}+E_{0111}+E_{0210}+E_{0211}+E_{1000}+E_{1100},\\\nonumber
f &=& 3 F_{0010}+3 F_{0011}+F_{0110}+F_{0111}+\frac{5 }{4}F_{0210}+\frac{5}{4}
   F_{0211}+6 F_{1000}+2 F_{1100},\\\nonumber
 h & =&   5 [E_{0001},F_{0001}]+10 [E_{0010},F_{0010}]+7 [E_{0100},F_{0100}]+4 [E_{1000},F_{1000}],
\end{eqnarray}
The following vectors form a complete set of maximum weight vectors of the irreducible $sl_2$-submodules. They are of eigenvalues 1,5,5,4,3,2,1, respectively, under $\ad_h$.
\begin{eqnarray}
L_2 &=&  \frac{20}{13} E_{0010}-\frac{28 }{13}E_{0011}-\frac{76
   }{13} E_{0110}-\frac{28}{13} E_{0111}+\frac{38}{13} E_{0210}+\frac{2
  }{13} E_{0211}+\frac{32 }{13} E_{1000}-\frac{88 }{13}E_{1100},\\\nonumber L_3&=&
   \frac{39}{20} E_{2431}, ~~~ L_4= \frac{39}{20} E_{2431}+\frac{9
   }{4}E_{2432}, \\\nonumber L_5&=& E_{2321}+E_{2421}, ~~ L_6=  2 E_{1221}+6
   E_{1321}+E_{2210}-5 E_{2211},\\\nonumber L_7&=& -4 E_{221}+E_{1110}-5
   E_{1111}-E_{1210}+5 E_{1211}, \\\nonumber L_8&=& \frac{2 }{5}E_{0010}+2
   E_{0011}-\frac{6}{5} E_{0110}-\frac{14
   }{5}E_{0111}-\frac{1}{5}E_{0210}+E_{0211}-\frac{4}{5} E_{1100}.
\end{eqnarray}
 Then setting
\begin{eqnarray}
K_1&=& F_{2432}, ~~ K_2=  \frac{15 }{13}F_{2431}-F_{2432}\\\nonumber
K_3 & =&  \frac{39
   }{20}F_{0010}-\frac{39 }{20}F_{0011}-\frac{39 }{8}F_{0110}-\frac{273
  }{40} F_{111}+\frac{39 }{10}F_{0211}-\frac{39 }{10}F_{1100}\\\nonumber K_4 & =&
   \frac{39}{20} F_{0010}+\frac{204 }{5}F_{0011}+3 F_{0110}-\frac{51
  }{5}F_{0111}+\frac{9 }{2}F_{0210}+\frac{129 }{10}F_{0211}+\frac{9}{4}
   F_{1000}+\frac{48}{5} F_{1100}
\end{eqnarray}
The vectors $\Y_i=L_i+K_i$ are basis of the opposite Cartan subalgebra $\h'$. The normailzed bilinear form is given by $\bil {g_1}{g_2}=\frac{1}{216} \Tr{(g_1.g_2)}$. Then one can check that $\bil {\Y_i}{\Y_j}=6 \delta_{ij}$. The  basis $\gamma_i\in \g^f$ such that $\bil {\gamma_i}{L_j}=\delta_{ij}$ are given by the formula
\begin{eqnarray}
\gamma_1 & = & f,\\\nonumber \gamma_2 & = &  \frac{1677
   }{1120}F_{0010}-\frac{1833 }{1120}F_{0011}-\frac{923
  }{1120} F_{0110}+\frac{247 }{1120}F_{0111}+\frac{403
   }{2240}F_{0210}+\frac{247 }{2240}F_{0211}+\frac{39
   }{35}F_{1000}-\frac{143}{140} F_{1100}, \\\nonumber \gamma_3 & = & \frac{15
  }{13} F_{2431}-F_{2432},~~ \gamma_4= F_{2432},\\\nonumber  \gamma_5 &=& \frac{27}{10} F_{2321}+\frac{9
  }{10} F_{2421},~~ \gamma_6=\frac{5 }{16}F_{1221}+\frac{5 }{16}F_{1321}-\frac{3
   }{8}F_{2210}-\frac{7}{8} F_{2211}, \\\nonumber \gamma_7 & = & -\frac{15 }{28}F_{0221}-\frac{27
   }{28}F_{1110}-\frac{9 }{4}F_{1111}+\frac{9 }{28}F_{1210}+\frac{3
   }{4}F_{1211}, \\\nonumber \gamma_8 & = & -\frac{405}{112} F_{0010}+\frac{135 }{16}F_{0011}+\frac{75
   }{112}F_{0110}-\frac{375 }{112}F_{0111}+\frac{45
   }{224}F_{0210}+\frac{15 }{32}F_{0211}-\frac{15}{14} F_{1100}.
\end{eqnarray}
We write elements of  Slodowy slice in the form $Q=L_1+ \sum_{i=1}^8 z_i\gamma_i$. The restriction $P_i^Q$ of the invariant polynomials $P_i$ of degree $\nu_i+1$  is obtained from  taking the trace of the matrix $Q^{\nu_i+1}$. We can take  $P_1^Q=z_1$. The expression correspond to the invariant of maximal degree $P_4^Q$ is omitted  since its very large. We give instead $\partial_{z_4}P_4^Q$.
\begin{eqnarray}
P_2^Q& =& 744192 z_1^3+\frac{44928}{7} z_2 z_1^2-\frac{518400}{7} z_8
   z_1^2-\frac{866970}{49} z_2^2 z_1+\frac{923400}{49} z_8^2 z_1-5760
   z_6 z_1\\\nonumber
   & - & \frac{1600560}{49} z_2 z_8 z_1+\frac{228002463}{137200}
   z_2^3-\frac{9871875}{686} z_8^3+\frac{150984}{49}
   z_7^2-\frac{37986975}{1372} z_2 z_8^2+\frac{165888}{13} z_3\\\nonumber &-& 3456
   z_4-\frac{6786}{7} z_2 z_6-\frac{45734949}{2744} z_2^2
   z_8+\frac{78300}{7} z_6 z_8, \end{eqnarray}

\begin{eqnarray}
P_3^Q &=& 40799232 z_1^4+958464 z_2 z_1^3-11059200 z_8 z_1^3-\frac{80016768}{35}
   z_2^2 z_1^2+\frac{24883200}{7} z_8^2 z_1^2-860160 z_6
   z_1^2\\\nonumber &- & \frac{31000320}{7} z_2 z_8 z_1^2+\frac{209079702}{1225} z_2^3
   z_1-\frac{89910000}{49} z_8^3 z_1+\frac{2287872}{7} z_7^2
   z_1-\frac{134573400}{49} z_2 z_8^2 z_1\\\nonumber &+ & \frac{24772608}{13}
   z_3
   z_1-516096 z_4 z_1-109824 z_2 z_6 z_1-\frac{84159972}{49} z_2^2
   z_8 z_1+1267200 z_6 z_8 z_1+\frac{9587156553}{686000}
   z_2^4\\\nonumber & -& \frac{29615625}{343} z_8^4-\frac{87267375}{343} z_2
   z_8^3+25920 z_6^2-\frac{112320}{49} z_2 z_7^2-\frac{29362905}{343}
   z_2^2 z_8^2\\\nonumber & +&\frac{621000}{7} z_6 z_8^2+207360 z_2 z_3  -149760 z_2
   z_4+\frac{534378}{35} z_2^2 z_6+311040 z_5 z_7+\frac{537489459}{6860} z_2^3
   z_8\\\nonumber &+ & \frac{1296000}{49} z_7^2 z_8-\frac{3456000}{13} z_3
   z_8+\frac{60840}{7} z_2 z_6 z_8,
\end{eqnarray}

\begin{eqnarray}
 \partial_{z_4} P_4^Q&=& -4505960448 z_1^3-\frac{18242205696}{7} z_2 z_1^2+\frac{1094860800}{7} z_8 z_1^2+\frac{2043055872}{245} z_2^2 z_1\\\nonumber &-&\frac{410572800}{7} z_8^2 z_1
  + 12165120 z_6
   z_1+\frac{5782233600}{49} z_2 z_8 z_1+\frac{20251269324 }{1225}z_2^3+\frac{801900000 }{343}z_8^3\\\nonumber &+&\frac{87588864}{49} z_7^2
   -\frac{1209265200}{49} z_2
   z_8^2-\frac{76972032}{13} z_3
  - 5308416 z_4
  \\\nonumber
  &+&\frac{41019264}{7} z_2 z_6+\frac{7000116552}{343} z_2^2 z_8+5702400 z_6 z_8.
\end{eqnarray}

Our special coordinates $(t_1,\ldots,t_8)$ are given by
\begin{eqnarray}
t_1& =&z_1,~t_2=-\frac{1}{149760}\partial_{z_4}P_3^Q-\frac{224}{65}z_1=z_2,~t_i=z_i, ~i=5,6,7,8\\\nonumber
t_3& =& -\frac{13}{331776000}\partial_{z_4}P_4^Q+\frac{13}{216000} P_3^Q=z_3+~\mathrm{nonlinear~terms},\\\nonumber
t_4&=&-\frac{1}{6912000}\partial_{z_4}P_4^Q-\frac{29}{432000} P_3^Q=z_4+~\mathrm{nonlinear~terms}.
\end{eqnarray}
Writing the restriction of the invariant polynomials in these coordinates, the space $N$ of common equilibrium points is defined as the zero set of the following polynomials
\begin{eqnarray}
\partial_{t_5} P_3^Q & =& 311040 t_7, \\\nonumber \partial_{t_6} P_3^Q &=&  -\frac{1478412 }{35}t_2^2-\frac{2779920 t_8 }{7}t_2+\frac{1458000}{7} t_8^2+51840
   t_6-622080 t_1 t_8, \\\nonumber
\partial_{t_7} P_3^Q & =& 311040 t_5-\frac{1866240}{7} t_1 t_7-\frac{9401184}{49} t_2 t_7+\frac{5598720}{49}  t_7
   t_8, \\\nonumber
\partial_{t_8} P_3^Q & =& 58844160 t_1^3+27248832 t_2 t_1^2-4147200 t_8 t_1^2-\frac{16116516}{35} t_2^2
   t_1+\frac{25758000}{7} t_8^2 t_1\\\nonumber &-& 622080 t_6 t_1+\frac{84240}{7} t_2 t_8
   t_1+\frac{31300659 }{980}t_2^3 -\frac{66825000 }{49}t_8^3+\frac{2799360
   }{49}t_7^2\\\nonumber &-&\frac{21718125}{49} t_2 t_8^2-\frac{3456000 }{13}t_3-\frac{2779920}{7} t_2
   t_6+\frac{38534535}{49} t_2^2 t_8+\frac{2916000}{7}  t_6 t_8. \\\nonumber
\end{eqnarray}
The local bihamiltonian structure is polynomial in  $t_1,t_2,t_3,t_4$ and $t_8$, where $t_8$ is a solution of a cubic equation.  The Potential of the Frobenius  structure  in  the flat coordinates $(s_1,s_2,s_3,s_4)$ is
\begin{eqnarray}
\mathbb F & =&T^2 (\frac{664832691 s_1^5}{43750}+\frac{393797781 s_2
   s_1^4}{8750}+\frac{117925163577 s_2^2 s_1^3}{2240000}+\frac{31524548679 s_2^3
   s_1^2}{1280000}\\\nonumber & & -\frac{177147 s_3 s_1^2}{1820}+\frac{1411599235293 s_2^4
   s_1}{286720000}-\frac{59049 s_2 s_3 s_1}{1120}+\frac{8090133251733
   s_2^5}{22937600000}-\frac{255879 s_2^2 s_3}{35840}) \\\nonumber & +& T (\frac{81990638748
   s_1^6}{546875}+\frac{157687224903 s_2 s_1^5}{546875}+\frac{252845042697 s_2^2
   s_1^4}{875000}+\frac{5680343128707 s_2^3 s_1^3}{28000000}\\\nonumber & &-\frac{405324 s_3
   s_1^3}{2275}+\frac{41422089388329 s_2^4 s_1^2}{448000000}-\frac{150903}{350} s_2 s_3
   s_1^2+\frac{349410443449509 s_2^5 s_1}{17920000000}\\\nonumber & &-\frac{2075463 s_2^2 s_3
   s_1}{5600}+\frac{118472583689109 s_2^6}{81920000000}+\frac{675
   s_3^2}{1183}-\frac{5051241 s_2^3 s_3}{89600}) \\\nonumber & +& \frac{2446443495072
   s_1^7}{13671875}+\frac{8512750428624 s_2 s_1^6}{13671875}+\frac{1593096854076 s_2^2
   s_1^5}{1953125}+\frac{87566456228121 s_2^3 s_1^4}{175000000}
   \\\nonumber & +& \frac{391896144 s_3
   s_1^4}{284375}+\frac{1357381494479907 s_2^4 s_1^3}{5600000000}-\frac{700488 s_2 s_3
   s_1^3}{21875}+\frac{21326967621723933 s_2^5 s_1^2}{224000000000}\\\nonumber & -& \frac{95335461 s_2^2
   s_3 s_1^2}{175000}+ \frac{87348137456366631 s_2^6 s_1}{3584000000000}+\frac{16}{169}
   s_3^2 s_1+\frac{1}{2} s_4^2 s_1-\frac{505028277 s_2^3 s_3
   s_1}{1400000}\\\nonumber & +&\frac{120333341133594693 s_2^7}{57344000000000}-\frac{7}{13} s_2
   s_3^2-\frac{10700732367 s_2^4 s_3}{89600000}+s_2 s_3 s_4
\end{eqnarray}
 where $T$ is a solution of the following cubic equation

 \begin{eqnarray}
0& =&  T^3-\frac{15552}{625} T s_1^2-\frac{4563 }{2500}T s_2^2-\frac{8424}{625} T s_1
   s_2-\frac{213504 }{15625}s_1^3-\frac{270231 }{62500}s_2^3\\\nonumber & & -\frac{444132 s_1
   }{15625}s_2^2-\frac{516672 s_1^2 }{15625}s_2+\frac{256}{2925} s_3.
 \end{eqnarray}
Then the quasihomgeneity condition reads
\be
{1\over 3} \partial_{s_1}\mathbb F+{1\over 3} \partial_{s_2}\mathbb F+ \partial_{s_3}\mathbb F+ \partial_{s_4} \mathbb F=(3-{2\over 3})\mathbb F.
\ee
\section{Conclusions and remarks}

Consider a nilpotent element not of semisimple type and the associated Drinfeld-Sokolov bihamiltonian structure. Then  the space  of common equilibrium points is still well defined and probably possesses  a local bihamiltonian structure which admits a dispersionless limit. However, examples show that its leading term does not define  a flat pencil of metrics.

It is known that for each conjugacy class in the Weyl group one can construct Drinfeld-Sokolov hierarchy \cite{gDSh1} and, under some restrictions,  an accompanied bihamiltonian structure \cite{gDSh2}. This bihamiltonian structure agrees with the one used  in this article if the conjugacy class is regular \cite{DelFeher}.

In the case of a regular primitive conjugacy classes, we obtain a new local algebraic bihamiltonian  structure on the space of common equilibrium points. Since it defines an exact Poisson pencil, its central invariants are constants \cite{FalLor}. It will be interesting to calculate them and  find if they are equal. In this case the bihamiltonian structure will be related to the topological hierarchy associated with the algebraic Frobenius structure \cite{DZ}.  This topological  hierarchy seems to be a reduction of the Drinfeld-Sokolov hierarchy (see \cite{mypaper4} for details on Dirac reduction of Hamiltonian equations).

 In future work, we will analyze the bihamiltonian structure associated to Drinfeld-Sokolov hierarchy for a primitive non-regular conjugacy class. Hoping, this will lead to algebraic Frobenius structure not covered in this article.

\noindent{\bf Acknowledgments.}

 The author thanks Boris Dubrovin for posting him this problem and for encouragement, support and useful discussions. The author also thanks Di Yang for stimulating discussions and anonymous reviewers whose comments/suggestions helped improve and clarify this article.   A part of this work was done during the author visits to the Abdus Salam International Centre for Theoretical Physics (ICTP) and the International School for Advanced Studies (SISSA) through the years 2014-2017. This work was also funded by the internal grant of Sultan Qaboos University (IG/SCI/DOMS/15/04).

~~~~~~~~
~~~~~~~~

\noindent Yassir Dinar

\noindent Department of Mathematics

\noindent College of Science

\noindent Sultan Qaboos University

\noindent dinar@squ.edu.om.

\end{document}